\documentclass[12pt]{article}

\usepackage[T1]{fontenc}
\usepackage[utf8]{inputenc}
\usepackage{amsmath,amssymb,epsfig,bbm}
\usepackage{stmaryrd,stix}
\usepackage[colorlinks]{hyperref}
\usepackage{comment}
\usepackage{color}
\usepackage{textcomp}
\usepackage[french]{babel}

\usepackage[acronym,toc]{glossaries}
\usepackage{glossary-longragged}
\usepackage{glossary-superragged}
\newglossary[slg]{symbols}{sym}{sbl}{List of Symbols}

\usepackage[textsize=small]{todonotes}
\usepackage{varwidth}

\usepackage{url}

\usepackage{enumitem}
\setlist[itemize]{labelindent=*, leftmargin=.5 truecm,nosep}

\usepackage[nottoc]{tocbibind}

\usepackage{amsthm}

\theoremstyle{plain}
\newtheorem{defi}{D\'efinition}
\newtheorem{prop}[defi]{Proposition}
\newtheorem{theo}[defi]{Th\'eor\`eme}
\newtheorem{conj}[defi]{Conjecture}
\newtheorem{lemm}[defi]{Lemme}
\newtheorem{coro}[defi]{Corollaire}

\theoremstyle{definition}
\newtheorem{rema}[defi]{Remarque}
\newtheorem{exem}[defi]{Exemple}
\newtheorem{exems}[defi]{Exemples}

\newcommand{\bdefi}{\begin{defi}}
\newcommand{\edefi}{\end{defi}}
\newcommand{\bprop}{\begin{prop}}
\newcommand{\eprop}{\end{prop}}
\newcommand{\btheo}{\begin{theo}}
\newcommand{\etheo}{\end{theo}}
\newcommand{\btheofr}{\begin{theofr}}
\newcommand{\etheofr}{\end{theofr}}
\newcommand{\blemm}{\begin{lemm}}
\newcommand{\brema}{\begin{rema}}
\newcommand{\erema}{\end{rema}}
\newcommand{\bexer}{\begin{exem}}
\newcommand{\eexer}{\end{exem}}
\newcommand{\bexem}{\begin{exem}}
\newcommand{\eexem}{\end{exem}}
\newcommand{\bexems}{\begin{exems}}
\newcommand{\eexems}{\end{exems}}
\newcommand{\bconj}{\begin{conj}}
\newcommand{\econj}{\end{conj}}
\newcommand{\elemm}{\end{lemm}}
\newcommand{\bcoro}{\begin{coro}}
\newcommand{\ecoro}{\end{coro}}
\newcommand{\dem}{\noindent{\bf Démonstration. }}

\usepackage{mathrsfs}
\renewcommand\mathcal{\mathscr}

\newcommand{\A}{{\cal A}}

\newcommand{\E}{{\cal E}}

\newcommand{\G}{{\cal G}}

\newcommand{\I}{{\cal I}}
\newcommand{\K}{{\cal K}}

\newcommand{\M}{{\cal M}}

\renewcommand{\O}{{\cal O}}
\renewcommand{\P}{{\cal P}}

\newcommand{\T}{{\cal T}}

\newcommand{\maths}[1]{{\mathbb #1}}  

\newcommand{\NN}{\maths{N}}

\newcommand{\QQ}{\maths{Q}}
\newcommand{\RR}{\maths{R}}

\newcommand{\ZZ}{\maths{Z}}

\newcommand{\fff}{{\mathfrak f}}
\renewcommand{\ggg}{{\mathfrak g}}

\newcommand{\OOO}{{\mathfrak O}}

\newcommand{\bs}{\backslash}

\newcommand{\Ga}{\Gamma}
\newcommand{\ov}[1]{{\overline{#1}}} 
\newcommand{\ra}{\rightarrow}

\newcommand{\wh}[1]{{\widehat{#1}}}

\newcommand{\card}{{\operatorname{Card}}}

\newcommand{\covol}{\operatorname{Covol}}

\newcommand{\cqfd}{\hfill$\Box$}

\newcommand{\gengeod}
{\operatorname{\widecheck{\G\,}\!\!}}

\newcommand{\Hom}{\operatorname{Hom}}

\newcommand{\id}{\operatorname{id}}

\newcommand{\rk}{\operatorname{rk}}

\newcommand{\ssm}{\!\smallsetminus\!}

\newcommand{\GL}{\operatorname{GL}}
\newcommand{\PGL}{\operatorname{PGL}}

\usepackage{interval}

\usepackage{mathtools}
\makeatletter
\DeclareRobustCommand\widecheck[1]{{\mathpalette\@widecheck{#1}}}
\def\@widecheck#1#2{%
    \setbox\z@\hbox{\m@th$#1#2$}%
    \setbox\tw@\hbox{\m@th$#1%
       \widehat{%
          \vrule\@width\z@\@height\ht\z@
          \vrule\@height\z@\@width\wd\z@}$}%
    \dp\tw@-\ht\z@
    \@tempdima\ht\z@ \advance\@tempdima2\ht\tw@ \divide\@tempdima\thr@@
    \setbox\tw@\hbox{%
       \raise\@tempdima\hbox{\scalebox{1}[-1]{\lower\@tempdima\box
\tw@}}}%
    {\ooalign{\box\tw@ \cr \box\z@}}}
\makeatother

\newcounter{fig}

\title{Minima successifs des réseaux et pentes des fibrés \\
  vectoriels sur les corps de fonctions globaux}
\author{Jean-Benoît Bost \and Fr\'ed\'eric Paulin}

\begin{document}
\bibliographystyle{alphanum}
\maketitle
\begin{abstract} Soient ${\bf C}$ une courbe projective lisse connexe
sur un corps fini, et $A$ l'algèbre affine de ses fonctions régulières
hors d'une place fixée de ${\bf C}$. Nous donnons des relations
précises entre les minima successifs de Mahler des $A$-réseaux normés
et les pentes de Harder-Narasimhan des fibrés vectoriels sur ${\bf C}$
via leur équivalence de catégorie.
\end{abstract}

\section{Introduction}
\label{sect:intro}

Soient $n\in\NN\ssm\{0\}$, ${\bf C}$ une courbe projective lisse
géométriquement connexe de genre $\ggg$ sur un corps fini $k_0$, $K$
son corps des fonctions, $v$ une place de $K$ de degré $\fff$, $\wh K$
le complété de $K$ en $v$, et $A$ l'algèbre affine de ${\bf C} \ssm
\{v\}$.

Un {\it $A$-réseau normé} est un couple $\overline{M}= (M,\|\;\|)$
d'un $A$-module $M$ projectif de type fini et d'une norme
ultramétrique $\|\;\|$ sur $M_{\wh K}=M\otimes_A\wh K$ (normalisée
pour être d'ensemble des valeurs exactement $\{0\}\cup
|k_0|^{\,\fff\,\ZZ}$ si $M_{\wh K}\neq \{0\}$).  Comme un analogue
projectif de l'équivalence de catégorie de Serre \cite[\S 48]{Serre55}
entre la catégorie des modules projectifs de type fini sur l'anneau de
coordonnées d'une variété algébrique affine $V$ et la catégorie des
faisceaux algébriques cohérents sur $V$, de nombreux travaux (de
Serre, Grothendieck, Atiyah, Harder, Narasimhan, Seshadri, etc) ont
établi une équivalence de catégorie $\overline{M}\mapsto
E^{\overline{M}}$ (préservant le rang) entre la catégorie des
$A$-réseaux normés et la catégorie des fibrés vectoriels (algébriques)
sur ${\bf C}$.

Notons $\lambda_1 (\,\overline{M}\,)\leq \lambda_2(\,\overline{M}\,)
\leq \dots\leq\lambda_n(\,\overline{M}\,)$ les minima successifs au
sens de Minkovski et Mahler d'un $A$-réseau normé $\overline{M}$ de
rang $n$. Nous renvoyons à la partie \ref{sect:resnorm} pour la
définition de ces minima successifs, voir \cite{Mahler41} lorsque
$\ggg=0$ et $\fff=1$, ainsi que \cite{Thunder94,RoyWal17,
  PoeRoy23}. Notons $\mu_1 (E)\geq \mu_2(E)\geq \dots\geq\mu_n(E)$ les
pentes de Harder-Narashiman d'un fibré vectoriel $E$ sur ${\bf C}$ de
rang $n$. Nous renvoyons à la partie \ref{sect:pentesHN} pour la
définition de ces pentes, due à \cite{HarNar75}, ainsi que par exemple
à \cite{Stuhler76, Grayson84, Bost96, Gaudron08, Andre09, Fargues10,
  Chen10a, Chen10b, Cornut18, Bost20, CheMor20, Grieve23, Li24}.

Le but de cette note est de revisiter le dictionnaire entre
$A$-réseaux normés et fibrés vectoriels (voir la partie
\ref{sect:fibvect}), et d'expliciter les relations précises entre les
minima successifs d'un $A$-réseau normé $\overline{M}$ et les pentes
de Harder-Narashiman de son fibré vectoriel associé
$E^{\overline{M}}$.  Dans le cadre des corps de nombres, ces relations
ont déjà de nombreuses applications en approximation diophantienne,
voir par exemple \cite{FalWus94,BreSax23, deSaxce23}. Pour tous les
entiers $p,q\in\ZZ$, nous noterons $\llbracket p,q\rrbracket=
[p,q]\cap\ZZ$ l'intervalle des entiers entre $p$ et $q$.

\btheo \label{theo:mainintro} Il existe des constantes explicites
$c_{\ggg,\fff,n}, c'_{\ggg,\fff,n}\geq 0$ telles que pour tout
$A$-réseau normé $\overline{M}$ de rang $n$, il existe des
sous-$A$-modules normés $\overline{M}_1, \dots, \overline{M}_n$ de
rang $1$ de $\overline{M}$ tels que nous ayons $M=M_1\oplus\dots\oplus
M_n$ et pour tout $i\in\llbracket1,n\rrbracket$,
\begin{equation}\label{eq:devissminsucc}
\big|\;\lambda_i(\,\overline{M}\,)-
\lambda_1(\,\overline{M}_i\,)\,\big|\leq c_{\ggg,\fff,n}\;,
\end{equation}
\begin{equation}\label{eq:relatpenminsuccintro}
\textstyle
\lambda_i(\,\overline{M}\,)\in \Big\llbracket\,
\big\lceil\textstyle{-\frac{1}{\fff}}\,\mu_i(E^{\overline{M}})\,\big\rceil
-\,c'_{\ggg,\fff,n},
\big\lfloor\textstyle{-\frac{1}{\fff}}\,\mu_i(E^{\overline{M}})\,\big\rfloor
+\,c'_{\ggg,\fff,n}\,\Big\rrbracket\;.
\end{equation}
\etheo

Nous renvoyons aux théorèmes \ref{theo:reducedbasis} et
\ref{theo:relatpenminsucc} pour des versions précisant les constantes
$c_{\ggg,\fff,n}$ et $c'_{\ggg,\fff,n}$. Celles-ci s'annulent lorsque
$\ggg=0$ et $\fff=1$, et les inégalités \eqref{eq:devissminsucc} et
\eqref{eq:relatpenminsuccintro} deviennent alors des égalités:
$\lambda_i(\,\overline{M}\,)= \lambda_1(\,\overline{M}_i\,)$ et
$\lambda_i(\,\overline{M}\,)= -\mu_i(E^{\overline{M}}\, )$.  Lorsque
${\bf C}$ est la droite projective sur $k_0$, tout fibré vectoriel sur
${\bf C}$ est somme directe de fibrés en droites (voir par exemple
\cite{Grothendieck57}). C'est en quantifiant de manière précise, dans
le théorème \ref{theo:quasiscind}, le défaut de véracité de cette
affirmation sur des courbes générales, puis en utilisant le
dictionnaire dans l'autre sens, que nous pourrons démontrer ce
théorème à la fin de cette note.

En prenant dans chaque $\overline{M}_i$ un élément $b_i$ non nul de
norme minimale, le théorème \ref{theo:mainintro} permet (voir la fin
de la partie \ref{sect:resnorm}) de construire une $A$-base
$(b_1,\dots, b_n)$ d'un sous-module libre de $\overline{M}$ d'indice
contrôlé et réalisant à constante contrôlée près les minima successifs
de Minskovski-Mahler de $M$. Ce résultat généralise la construction
d'une $A$-base de Mahler de $\overline{M}$ lorsque $\ggg=0$ et
$\fff=1$, voir \cite{Mahler41} et \cite[\S 1]{Lenstra85} pour une
version algorithmique. Dans l'analogie entre corps de fonctions et
corps de nombres, le théorème \ref{theo:mainintro} est une version
effective du pendant pour les corps de fonctions de la théorie de la
réduction des réseaux euclidiens (dont la construction de bases
réduites de Minkovski ou de Korkin-Zolotarev, voir par exemple
\cite{LagLenSch90}). Dans l'interprétation par les normes
ultramétriques des quotients d'immeubles de Bruhat-Tits sur les
$K$-groupes réductifs (voir \cite{Harder69} pour l'aspect des pentes),
ce résultat s'insère donc dans l'analogie avec la théorie de la
réduction des quotients arithmétiques associés aux $\QQ$-groupes
réductifs (voir \cite{Stuhler76, Grayson84} pour l'aspect des
pentes). Dans ce dernier cadre, les constructions qui sous-tendent les
opérations de troncature mises en oeuvre dans la formule des traces
d'Arthur-Selberg (voir \cite[\S 6] {Arthur78}) apparaissent comme une
généralisation du formalisme des pentes des réseaux euclidiens, comme
remarqué dans \cite[\S 5.1]{Bost20}.

Nous comprenons que cette note pourrait être généralisée au cas des
complétions $S$-adiques pour $S$ un ensemble fini de places de $K$
(avec applications en géométrie des nombres $S$-adiques comme
développé dans \cite{KleShiTom17}) ou de manière plus générale aux
courbes adéliques (comme développé dans \cite{CheMor20}). Mais nous
avons choisi de rester dans le cadre d'une seule place afin de clarifier
l'exposition en vue d'applications dans \cite{BanKimLimPau24} en
géométrie des nombres paramétrique sur des corps de fonctions globaux
généraux, généralisant ceux de Roy-Waldschmidt \cite{RoyWal17} lorsque
$\ggg=0$ et $\fff=1$ (voir aussi \cite{PoeRoy23}).

\medskip
\noindent{\small {\it Remerciements :} Nous remercions
  Emanuele Macri pour la démonstration du lemme \ref{lem:macri}. }

\section{Les minima successifs de Minkovski des réseaux}
\label{sect:resnorm}

\noindent{\bf Notations sur les corps globaux de fonctions. } Nous
renvoyons par exemple à \cite{Rosen02,Goss96} pour ces rappels. Soient
$k_0$ un corps fini et $q_0$ son ordre. Soient ${\bf C}$ une courbe
projective lisse géométriquement connexe sur $k_0$, et $\ggg$ son
genre. Soit $K=k_0({\bf C})$ le corps des fonctions de ${\bf C}$.
Nous fixons un point fermé $P_\infty$ de ${\bf C}$, appelé le {\it
  point à l'infini} de ${\bf C}$. Notons $v$ la valuation (discrète
normalisée) associée à $P_\infty$, $\fff=\deg(P_\infty)=
[k_0(P_\infty):k_0]$ le degré de $P_\infty$, $\OOO$ l'anneau local en
$P_\infty$ dans $K$, $\pi$ une uniformisante dans $\OOO$ de sorte que
$v(\pi)=1$, $k= \OOO/\pi\OOO$ le corps résiduel de $\OOO$, $q=
q_0^{\;\fff}$ l'ordre de $k$, $|\cdot|= q^{\;-v(\cdot)}$ la valeur
absolue normalisée associée à $P_\infty$, $\wh K$ le complété de $K$
pour $|\cdot|$, $\wh \OOO$ celui de $\OOO$, et $A$ l'algèbre affine de
la courbe affine ${\bf C}^\circ={\bf C} \ssm\{P_\infty\}$. Rappelons
que $A$ est un anneau de Dedekind.

\medskip
\noindent{\bf Rappels d'algèbre linéaire normée ultramétrique. } Soit
$n\in\NN\ssm\{0\}$.  Soit $V$ un espace vectoriel de dimension $n$ sur
le corps valué $\wh K$. Soit $\|\cdot\|$ une norme (ultramétrique) sur
$V$, dont nous noterons $B_{\|\cdot\|}(x,r)$ les boules fermées de
centre $x$ et de rayon $r$.  Nous dirons que $\|\cdot\|$ est {\it
  entière} si l'ensemble de ses valeurs est $\{0\}\cup q^\ZZ$. Le
groupe linéaire $\GL(V)$ agit (à gauche) sur l'ensemble
$\operatorname{Norm}_\ZZ(V)$ des normes entières de $V$, par
précomposition par l'inverse $(g,\|\cdot\|)\mapsto g\|\cdot\|:x\mapsto
\|g^{-1}x\|$. Si $n\geq 2$, dans le modèle des classes d'homothéties
des normes sur $\wh K^n$ de l'immeuble affine $\I(\PGL_n,\wh K)$ de
Bruhat-Tits de $\PGL_n$ sur $\wh K$ (voir par exemple
\cite{GolIwa63,Parreau00}), l'ensemble $\operatorname{Norm}_\ZZ(\wh
K^n)$ est exactement l'ensemble des sommets de $\I(\PGL_n,\wh K)$ :la
condition d'être entière pour une norme en fait un représentant
canonique de sa classe d'homothétie.

En suivant \cite[\S II.1]{Weil95}, nous dirons qu'une $\wh K$-base
$(e_1,\dots, e_n)$ de $V$ est {\it orthonormée} pour $\|\cdot\|$ si
pour tous les $x_1,\dots, x_n\in\wh K$, nous avons $\|x_1e_1+\dots+
x_ne_n\|=\max_{1\leq i\leq n}|x_i|$. Il découle par exemple de
\cite[Prop.~1.1]{GolIwa63} qu'une norme sur $V$ admet une $\wh K$-base
orthonormée si et seulement si elle est entière. Notons encore
$\|\cdot\|$ la norme duale de $\|\cdot\|$ sur l'espace vectoriel dual
$V^*$ de $V$ définie par $\|\ell\|=\sup_{x\in V\ssm \{0\}}
\frac{|\ell(x)|}{\|x\|}$. Pour tout $j\in\NN$, nous noterons encore
$\|\cdot\|$ la norme sur la $j$ème puissance extérieure $\bigwedge^jV$
de $V$ définie, pour tout $w\in \bigwedge^j V$, par
\[
\|w\|=\sup_{\ell_1,\,\dots,\,\ell_j\,\in\, V^*\;\ssm\;\{0\}}
\frac{|(\ell_1\wedge\dots\wedge \ell_j)(w)|}
{\|\ell_1\|\dots\|\ell_1\|}\;.
\]
Si la norme $\|\cdot\|$ de $V$ est entière, alors pour toute $\wh
K$-base orthonormée $(e_1,\dots, e_n)$ de $(V,\|\cdot\|)$, la $\wh
K$-base $(e_{i_1}\wedge \dots \wedge e_{i_j})_{1\leq i_1<\dots<i_j\leq
  n}$ de $\bigwedge^jV$ est orthonormée pour la norme $\|\cdot\|$ de
$\bigwedge^jV$.

Rappelons qu'un {\it $\wh \OOO$-réseau} $L$ de $V$ est un $\wh
\OOO$-sous-module de type fini (donc libre car $\wh \OOO$ est
principal) engendrant $V$ sur $\wh K$ (ou, de manière équivalente,
étant compact-ouvert). Le groupe linéaire $\GL(V)$ agit transitivement
sur l'ensemble $\operatorname{Res}_{\wh \OOO}(V)$ des $\wh
\OOO$-réseaux de $V$. Si $V={\wh K}^n$, le stabilisateur dans
$\GL(V)=\GL_n(\wh K)$ du $\wh \OOO$-réseau produit ${\wh \OOO}^n$ est
$\GL_n(\wh \OOO)$. Nous notons $L_{\|\cdot\|}$ la boule unité fermée
${\ov B}_{\|\cdot\|}(0,1)$ de $\|\cdot\|$, qui est un $\wh
\OOO$-réseau de $V$. Nous notons alors $m_{\|\cdot\|}$ la mesure de
Haar du groupe additif localement compact $V$ normalisée par
$m_{\|\cdot\|} (L_{\|\cdot\|}) =1$.  L'application $\|\cdot\|\mapsto
L_{\|\cdot\|}$ est une bijection $\GL(V)$-équivariante de
$\operatorname{Norm}_\ZZ (V)$ dans $\operatorname{Res}_{\wh \OOO}(V)$
d'inverse l'application
\[
L\mapsto \big(x\mapsto \|x\|_L=
\inf\{|\lambda|:\lambda\in\wh K, x\in \lambda L\}\big)
\]
qui à un $\wh \OOO$-réseau $L$ de $V$ associe sa jauge de Minkovski
$\|\cdot\|_L$ (voir par exemple \cite[\S 1.1.7]{CheMor20},
\cite[Prop.~1.5.2]{BosCha24}).

Rappelons qu'un {\it $A$-réseau} $\Lambda$ de $V$ est un
$A$-sous-module de $V$ discret de type fini (donc projectif car sans
torsion) qui engendre $V$ sur $\wh K$. Nous avons une identification
canonique $\Lambda\otimes_A K=V$. L'espace topologique quotient
$V/\Lambda$ est compact. Nous le munissons de la mesure induite par
$m_{\|\cdot\|}$ (telle que la projection canonique $V\ra V/\Lambda$
préserve localement la mesure), encore notée $m_{\|\cdot\|}$. Nous
appellerons {\it covolume} de $\Lambda$ pour la norme $\|\cdot\|$ la
masse totale $\covol_{\|\cdot\|}(\Lambda) =m_{\|\cdot\|}(V/\Lambda)$
de $V/\Lambda$.  Le $A$-réseau $A^n$ de ${\wh K}^n$, muni de la norme
rendant orthonormée la ${\wh K}$-base canonique de ${\wh K}^n$, est de
covolume égal à
\begin{equation}\label{eq:covolRv}
  \covol(A^n)= q_0^{\;n(\ggg-1)}
\end{equation}
(voir par exemple \cite[Lem. 14.4]{BroParPau19}). Nous appellerons {\it
  covolume normalisé} d'un $A$-réseau $\Lambda$ de $V$ pour la norme
$\|\cdot\|$ la quantité
\[
\overline{\covol}_{\|\cdot\|}(\Lambda)
=\frac{\covol_{\|\cdot\|} (\Lambda)}{\covol(A^n)}\,.
\]
Nous dirons que $\Lambda$ est {\it unimodulaire} pour la norme
$\|\cdot\|$ si $\overline{\covol}_{\|\cdot\|}(\Lambda)=1$. Notons que
pour tout $g\in \GL(V)$, nous avons
\[
\covol_{\|\cdot\|}(g\Lambda)=|\det g|\; \covol_{\|\cdot\|}
(\Lambda)\,.
\]
En particulier, si la norme $\|\cdot\|$ est entière, si $\Lambda$ est
libre sur $A$, et si $(b_1,\dots,b_n)$ est une $A$-base de $\Lambda$,
alors
\begin{equation}\label{eq:covolnormwedge}
  \overline{\covol}_{\|\cdot\|}\,(\Lambda)=
  \|\,b_1\wedge\dots\wedge b_n\,\|\;.
\end{equation}

\medskip
\noindent{\bf Minimaux successif des $A$-réseaux normés. } Soient
$n\in\NN$, $V$ un espace vectoriel sur $\wh K$ de dimension $n$ muni
d'une norme $\|\;\|$, et $\Lambda$ un $A$-réseau de $V$. Pour tout
$i\in \llbracket 1,n\rrbracket$, le {\it $i$-ème minimum} de $\Lambda$
pour la norme $\|\;\|$ est
\[
\lambda_i (\Lambda,\|\;\|)=\min\{\rho\in\RR : \dim_{\wh K}
(\operatorname{vect}_{\wh K}(B_{\|\;\|}(0,q^{\rho})\cap \Lambda))\geq i\}.
\]
Cette définition est en fait l'image par $\log_{q}$ de la définition
la plus fréquement utilisée (voir par exemple \cite{Mahler41}).
Cette convention simplifiera les liens avec les pentes, voir aussi
\cite{Thunder94} qui utilise la même convention (utile en particulier
pour considérer des corps de constante infinis). Le premier minimum
$\lambda_1(\Lambda,\|\;\|)$ est aussi appelé la {\it systole} de
$(\Lambda,\|\;\|)$.

Le minimum ci-dessus existe par la compacité des boules de $V$ et par
la discrétude des $A$-réseaux. Nous avons $\lambda_i(\Lambda,\|\;\|)
\in \ZZ$ si $\|\;\|$ est entière. Nous avons clairement
\[
\lambda_1
(\Lambda,\|\;\|)\leq \dots \leq \lambda_n(\Lambda,\|\;\|)\,.
\]
Pour tous les $i\in \llbracket 1,n\rrbracket$ et $g\in\GL(V)$, nous
avons $\lambda_i (g\Lambda, g\|\;\|) =\lambda_i (\Lambda,\|\;\|)$.
Puisque $\GL(V)$ agit transitivement sur $\operatorname{Norm}_\ZZ(V)$
ainsi que sur $\operatorname{Res}_{\wh \OOO}(V)$, cette propriété
d'invariance permet d'ou bien fixer le $A$-réseau $\Lambda$ et de
varier la norme entière $\|\;\|$, ou bien le contraire.

Si $\|\;\|'$ est une autre norme sur $V$, et si $c,c'\in\RR$ vérifient
$q^{c}\|x\|\leq \|x\|'\leq q^{c'}\|x\|$ pour tout $x\in V$, alors pour
tout $i\in \llbracket 1,n\rrbracket$, nous avons
\begin{equation}\label{eq:comparnormsuccmin}
\lambda_i (\Lambda,\|\;\|)+c\leq \lambda_i (\Lambda,\|\;\|')\leq
\lambda_i (\Lambda,\|\;\|)+c'\;.
\end{equation}

Les inégalités suivantes, qui découlent de
\cite[Theo.~4.4]{KleShiTom17} et de la formule \eqref{eq:covolRv},
voir aussi \cite[page 489]{Mahler41} lorsque $\ggg=0$ et $\fff=1$,
\cite[Theo.~2.1]{KimLimPau23} ainsi que \cite[Theo 4.1]{PoeRoy23} et
ses références, sont une version pour les corps de fonctions du
théorème de Minkowski sur les corps convexes des espaces euclidiens~:
si $\|\;\|$ est entière (voir \cite{BanKimLimPau24} pour l'ajustement
des constantes dans le cas des normes quelconques), alors
%
\begin{equation}\label{eq:Minkovski2}
\log_{q}\overline{\covol}_{\|\;\|}(\Lambda)\leq
\sum_{i=1}^n\lambda_i(\Lambda,\|\;\|)\leq \log_{q}
\covol_{\|\;\|}(\Lambda)+n\;.
\end{equation}

Une des motivations principales de cette note, qui sera utilisée dans
\cite{BanKimLimPau24}, est le résultat suivant, qui donne l'existence
dans tout $A$-réseau d'un sous-$A$-réseau libre muni d'une $A$-base
bien adaptée aux minimaux successifs. Il est dû à Mahler \cite[page
  489]{Mahler41} lorsque $\ggg=0$ et $\fff=1$. Nous le démontrerons à
la fin de cette note.

\btheo\label{theo:reducedbasis} Si la norme $\|\;\|$ de $V$ est
entière, pour tout $A$-réseau $\Lambda$ de $V$, il existe une
décomposition en somme directe de $A$-modules $\Lambda= \Lambda_1
\oplus\dots \oplus\Lambda_n$, où $\Lambda_1, \dots, \Lambda_n$ sont
des sous-$A$-modules de rang $1$ de $\Lambda$ tels que nous ayons
$\lambda_1(\Lambda_1,\|\;\|_{\mid \wh K\Lambda_1})\leq \dots\leq
\lambda_1(\Lambda_n,\|\;\|_{\mid \wh K\Lambda_n})$ et, pour tout
$i\in\llbracket1,n\rrbracket$,
\[
\big|\;\lambda_i(\Lambda,\|\;\|)-
\lambda_1(\Lambda_i,\|\;\|_{\mid \wh K\Lambda_i})\;\big|\leq c_{\ggg,\fff,n}=
\Big\lfloor\,\frac{n(n-1)}{2\,\fff}(4\ggg+3\fff-3)\,\Big\rfloor\;.
\]
\etheo

En notant $b_1\in \Lambda_1,\dots, b_n\in \Lambda_n$ des éléments tels
que $\|b_i\|=\lambda_1(\Lambda_i,\|\;\|_{\mid \wh K\Lambda_i})$, alors
$\Lambda_0=A b_1+\dots +A b_n$ est un sous-$A$-réseau libre de
$\Lambda$ (de $A$-base $(b_1,\dots, b_n)$), dont l'indice
$[\Lambda:\Lambda_0]= \frac{\covol_{\|\;\|}(\Lambda_0)}
{\covol_{\|\;\|} (\Lambda)}$ est uniformément borné (par une constante
ne dépendant que de $\ggg, \fff, n$) par la formule
\eqref{eq:Minkovski2}. Nous renvoyons à \cite{BanKimLimPau24} pour
d'autres précisions sur ce théorème et pour des applications en
géométrie des nombres paramétrique sur des corps de fonctions
généraux, la $A$-base $(b_1,\dots, b_n)$ de $\Lambda_0$ ci-dessus
remplaçant la $A$-base de Mahler de $\Lambda$ utilisée dans
\cite{RoyWal17} lorsque $\ggg=0$ et $\fff=1$.

\section{Les liens entre réseaux normés et fibrés vectoriels}
\label{sect:fibvect}

Nous renvoyons par exemple à \cite{Serre55} et \cite[\S
  II.2.1]{Serre83} pour cette partie. Soit $n\in\NN$.

\medskip\noindent{\bf $A$-réseaux normés. }
Nous appelons {\it $A$-réseau normé} de rang $n$ tout couple
$\overline{M}= (M,\|\cdot\|)$ constitué d'un $A$-module $M$ projectif
de rang $n$ et d'une norme entière $\|\cdot\|$ sur l'espace vectoriel
$M_{\wh K}=M\otimes_A\wh K$ de dimension $n$ sur le corps valué
localement compact $\wh K$.  Remarquons que $M$ est un $A$-réseau de
$M_{\wh K}$ au sens de la partie \ref{sect:resnorm}. Nous noterons
$\overline{\covol}\, (\,\overline{M}\,)$ le covolume normalisé du
$A$-réseau $M$ de $M_{\wh K}$ pour la norme $\|\cdot\|$ sur $M_{\wh
  K}$. Si $\overline{M}{}'= (M',\|\cdot\|')$ est un $A$-réseau normé
de rang $n$, un {\it morphisme} de $A$-réseaux normés $\phi:
\overline{M}\ra \overline{M}{}'$ est un morphisme de $A$-modules
$\phi:M\ra M'$ dont l'extension $\wh K$-linéaire $\phi_{\wh K}:M_{\wh
  K}\ra {M'}_{\wh K}$ est de norme d'opérateur $\interleave\,
\phi_{\wh K}\,\interleave =\sup_{x\in \wh{M}\;\ssm\;\{0\}}
\frac{\|\,\phi_{\wh K}(x)\,\|'}{\|x\|}$ au plus $1$.  Cette norme
d'opérateur sur le $\wh K$-espace vectoriel des applications linéaires
de $M_{\wh K}$ dans ${M'}_{\wh K}$ est aussi ultramétrique. Nous
noterons $\Hom(\,\overline{M}, \overline{M}{}'\,)$ le $k_0$-espace
vectoriel des morphismes de $A$-réseaux normés de $\overline{M}$ dans
$\overline{M}{}'$.

\medskip\noindent{\bf Fibrés vectoriels sur ${\bf C}$. } Notons
$\O_{\bf C}$ le faisceau structurel de la courbe projective ${\bf C}$
sur $k_0$. Un {\it fibré vectoriel} (algébrique) $E$ sur ${\bf C}$ est
un faisceau cohérent localement libre de $\O_{\bf C}$-modules sur
${\bf C}$. Un {\it sous-fibré vectoriel} de $E$ est un sous-faisceau
cohérent de $E$; c'est un fibré vectoriel sur ${\bf C}$. Le {\it
  saturé} $F^{\rm sat}$ d'un sous-fibré vectoriel $F$ de $E$ est le
noyau du morphisme de faisceaux quotient composé $E\ra (E/F)/M$ où $M$
est le sous-$\O_{\bf C}$-module de torsion du faisceau quotient
$E/F$. Un sous-fibré vectoriel $F$ de $E$ est {\it saturé} s'il est
égal à son saturé, c'est-à-dire si $E/F$ est un $\O_{\bf C}$-module
sans torsion, et alors $E/F$ est un fibré vectoriel sur ${\bf C}$.
L'intersection de deux sous-fibrés vectoriels saturés de $E$ est
saturée, mais la somme ne l'est pas forcément. La fibre $E_K$ de $E$
au point générique de ${\bf C}$ est un $K$-espace vectoriel de
dimension finie. Si $F$ est un sous-fibré vectoriel de $E$, alors
$F_K$ est un sous-espace vectoriel de $E_K$, et $F_K=(F^{\rm sat})_K$.
L'application $F\mapsto F_K$ est une bijection de l'ensemble des
sous-fibrés vectoriels saturés de $E$ dans l'ensemble des sous-espaces
vectoriels de $E_K$.

Soient $E$ et $F$ deux fibrés vectoriels sur ${\bf C}$. Un morphisme
de fibrés vectoriels de $E$ dans $F$ est un morphisme de faisceaux de
$E$ dans $F$ (son faisceau image n'est pas forcément saturé). Nous
définissons les opérations suivantes sur les fibrés vectoriels,
correspondantes par $E\mapsto E_K$ aux opérations éponymes sur les
$K$-espaces vectoriels de dimension finie : le dual
$E^{\,\widecheck{\;}}$, la $k$-ème puissance tensorielle $E^{\otimes
  k}$ pour tout $k\in\NN$ (et $E^{\otimes -k}=
(E^{\,\widecheck{\;}})^{\otimes k}$), la $k$-ème puissance extérieure
$\Lambda^kE$ de $E$, la somme directe $E\oplus F$, le produit
tensoriel $E\otimes F$ et le fibré des morphismes $\Hom_{\O_{\bf C}}
(E,F)$. Le {\it rang} du fibré vectoriel $E$ est
\[
\rk E =\dim_K E_K\,.
\]
Nous avons $\rk(E\oplus F)=\rk E + \rk F$, $\rk(E\otimes F)=(\rk
E)(\rk F)$ et, pour tout $k\in\ZZ$, $\rk E^{\otimes k}=(\rk E)^{|k|}$.
Un {\it fibré en droites} sur ${\bf C}$ est un fibré vectoriel sur
${\bf C}$ de rang $1$. Si $n$ est le rang de $E$, nous notons $\det E
= \bigwedge^nE$ le fibré en droites {\it déterminant} de $E$. Si
$D=\sum_{i\in I} n_iP_i$ est le diviseur de n'importe quelle section
rationnelle non nulle de $\det E$, alors le {\it degré} de $E$ est le
degré de $D$, donc
\[
\deg E  = \deg(\det E)= \deg D=\sum_{i\in I} n_i\,[k_0(P_i) :k_0]\,.
\]
Si $F$ est un sous-fibré vectoriel de $E$, alors $\rk F=\rk F^{\rm sat}$,
\begin{equation}\label{eq:satcrois}
\deg F\leq \deg F^{\rm sat}
\end{equation}
avec égalité si et seulement si $F= F^{\rm sat}$, et si $F$ est
saturé, nous avons 
\[
\deg E=\deg F+\deg (E/F)\,.
\]
Lorsque $E$ et $F$ sont des fibrés en droites sur ${\bf C}$, nous
avons $\deg (E\otimes F) = \deg E+\deg F$ et $\deg
E^{\,\widecheck{\;}} = -\deg E$; pour tout $k\in\ZZ$, nous avons alors
\[
\deg E^{\otimes k}=k\deg E\,.
\]
Si $m$ est le degré de $F$, nous avons $\Lambda^{nm}(E\otimes F)=
(\Lambda^{n}E)^{\otimes m} \otimes (\Lambda^{m}F)^{\otimes n}$, donc
\begin{equation}\label{eq:degprodtens}
  \deg(E\otimes F)= (\rk F)\deg(E) +(\rk E) \deg(F)\,.
\end{equation}

Nous utiliserons plusieurs fois le lemme classique suivant. Notons
$\P$ l'ensemble des points fermés de ${\bf C}$. Soit $\T$ un module de
type fini de torsion sur un anneau de valuation discrète $\A$ d'idéal
maximal $\M$. Par le théorème des diviseurs élémentaires, il existe
$n\in\NN$ et $a_1,\dots, a_n\in\NN\ssm\{0\}$ tel que $\T$ soit
isomorphe à $\prod_{i=1}^n\A/\M^{a_i}$. Par définition, la {\it
  longueur} $\operatorname{lg}\,(\T\,)$ de $\T$ est $\operatorname{lg}
\,(\T\,)=\sum_{i=1}^n a_i$.

\blemm\label{lem:macri}
Soit $E'$ un sous-fibré vectoriel d'un fibré vectoriel $E$
sur ${\bf C}$. Si $E$ et $E'$ ont le même rang, alors
\[
\deg E-\deg E'=\sum_{P\in\P}[k_0(P):k_0]\operatorname{lg}\,(E_P/E'_P)\,.
\]
\elemm

\dem Rappelons que pour tout $P\in\P$, la fibre $\O_{{\bf C},P}$ de
$\O_{{\bf C}}$ en $P$ est un anneau de valuation discrète, dont nous
noterons $\M_{{\bf C},P}$ l'idéal maximal, et que les fibres $E_{P}$
de $E$ et $E'_{P}$ de $E'$ en $P$ sont des $\O_{{\bf C},P}$-modules
sans torsion (libre). Notons $T$ le faisceau cohérent quotient $E/E'$,
qui est de support de dimension $0$ car $\rk E =\rk E'$. Pour tout
$P\in\P$, notons $T_P=E_P/E'_P$ la localisation de $T$ en $P$, de
sorte que $T=\oplus_{P\in\P}T_P$. Pour tout faisceau cohérent $F$ sur
${\bf C}$, notons $\chi({\bf C},F)=\sum_{i\in\NN}
(-1)^i\dim_{k_0}H^i({\bf C},F)$ la caractéristique d'Euler de $F$. Si
$F$ est un fibré vectoriel, alors $\deg F=\chi({\bf C},F)-(\rk
F)\chi({\bf C},\O_{\bf C})$, voir par exemple \cite[\S
  1.2]{HuyLeh10}. Puisque $E$ et $E'$ ont le même rang, par
l'additivité de la caractéristique d'Euler et puisque le support de
$T$ est de dimension $0$, nous avons
\begin{align}
  \deg E-\deg E'&=\chi({\bf C},E)-\chi({\bf C},E')=\chi({\bf C},T)=
  \dim_{k_0}H^0({\bf C},T)\nonumber\\&=\sum_{P\in\P}\dim_{k_0}T_P=
  \sum_{P\in\P}[k_0(P):k_0]\dim_{k_0(P)}T_P\,.\label{eq:macri}
\end{align}
Si $n_P\in\NN$ et $a_1,\dots, a_{n_P}\in\NN\ssm\{0\}$ sont tels que
$T_P$ soit isomorphe à $\prod_{i=1}^{n_P}\O_{{\bf C},P}/\M_{{\bf C},
  P}^{a_i}$, alors $\dim_{k_0(P)}T_P= \sum_{i=1}^{n_P} a_i=
\operatorname{lg}\,(T_P)$. Le résultat en découle.  \cqfd

\medskip\noindent{\bf Le dictionnaire entre $A$-modules projectifs
  normés et fibrés vectoriels sur ${\bf C}$. } Soit $\overline{M}$ un
$A$-réseau normé de rang $n$. Notons $E= E^{\overline{M}}$ le fibré
vectoriel de rang $n$ sur ${\bf C}$ qui est le sous-faisceau cohérent
localement libre du faisceau constant $M_{K}=M\otimes_A K$ sur ${\bf
  C}$, dont la fibre $E_{P_\infty}$ au point à l'infini $P_\infty$ est
le $\OOO$-réseau de $M_K$ égal à l'intersection avec $M_K$ du $\wh
\OOO$-réseau $L_{\|\cdot\|}= B_{\|\cdot\|} (0,1)$ de $M_{\wh K}$, et
dont l'espace des sections affines est $\Ga({\bf C}^\circ,E) =M$. Nous
avons alors
\begin{equation}\label{eq:relatsecttotal}
  \Ga({\bf C},E^{\overline{M}}) =
  \{s\in M:\|s\|\leq 1\}=M\cap L_{\|\cdot\|}\;.
\end{equation}
Réciproquement, pour tout fibré vectoriel $E'$ de rang $n$ sur
${\bf C}$, notons $\overline{M}{}^{E'}= (M^{E'}, \|\cdot\|^{E'})$ le
$A$-réseau normé de rang $n$ où $M^{E'}$ est le $A$-module projectif
$\Ga({\bf C}^\circ,E')$ et la norme $\|\cdot\|^{E'}$ est déterminée
par le fait que son $\wh \OOO$-réseau associé $L_{\|\cdot\|^{E'}}$
soit le $\wh \OOO$-réseau de $(M^{E'})_{\wh K}$ complétion (voir par
exemple \cite[\S 1.5.2]{BosCha24}) du $\OOO$-réseau $(E')_{P_\infty}$.

Soient $F'$ un fibré vectoriel de rang $n$ sur ${\bf C}$ et $\varphi:
E'\ra F'$ un morphisme de fibrés vectoriels entre $E'$ et $F'$. Notons
$\phi=\phi(\varphi): M^{E'}=\Ga({\bf C}^\circ,E')\ra M^{F'}= \Ga({\bf
  C}^\circ,F')$ le morphisme de $A$-modules $s\mapsto \varphi \circ s$
de post-composition par $\varphi$ des sections. Il vérifie que
l'application $\phi_{\wh K}:(M^{E'})_{\wh K}\ra (M^{F'})_{\wh K}$
envoie la boule unité pour $\|\cdot\|^{E'}$ dans la boule unité pour
$\|\cdot\|^{F'}$, donc est de norme d'opérateur $\interleave\,
\phi_{\wh K} \, \interleave$ au plus $1$.

Comme expliqué par exemple dans \cite[\S II.2.1]{Serre83}, le foncteur
$E'\mapsto \overline{M}{}^{E'}$ et $\varphi\mapsto \phi(\varphi)$ est
une équivalence de catégories d'inverse $\overline{M}\mapsto
E^{\overline{M}}$, qui vérifie les propriétés suivantes pour tous les
$A$-réseaux normés $\overline{M}=(M,\|\;\|)$ et
$\overline{M}{}'=(M',\|\;\|')$ de rang $n$.
\begin{enumerate}
\item\label{item1:dictionnaire} L'application de $\Hom_{\O_{\bf C}}
  (E^{\overline{M}}, E^{\overline{M}{}'})$ dans $\Hom(\,\overline{M},
  \overline{M}{}'\,)$ définie par $\varphi\mapsto \phi(\varphi)$ est
  un isomorphisme de $k_0$-espaces vectoriels de dimension finie, et
  les fibrés vectoriels $E^{\overline{M}}$ et $E^{\overline{M}{}'}$
  sur ${\bf C}$ sont isomorphes si et seulement si $\overline{M}$ et
  $\overline{M}{}'$ sont isomorphes.
\item\label{item2:dictionnaire}
Pour tout $\varphi\in\Hom_{\O_{\bf C}}(E^{\overline{M}},
E^{\overline{M}{}'})$, son image $\varphi(E^{\overline{M}})$ est un
sous-fibré vectoriel saturé de $E^{\overline{M}{}'}$ si et seulement si

$\bullet$~ d'une part l'image du morphisme $\phi(\varphi):M\ra M'$ est
saturée dans $M'$ (c'est-à-dire telle que $M'/\phi(\varphi)(M)$ est un
$A$-module sans torsion, ou de manière équivalente puisque $A$ est un
anneau de Dedekind, tel que $\phi(\varphi)(M)$ soit un facteur direct
de $M'$), et

$\bullet$~ d'autre part la norme image de $\|\;\|$ par $\phi_{\wh K}$,
définie par $x'\mapsto \min_{x\in \phi_{\wh K}^{-1}(x')} \|x\|$ pour
tout $x'\in \phi_{\wh K}(M_{\wh K})$, coïncide avec la restriction à
$\phi_{\wh K}(M_{\wh K})$ de la norme $\|\;\|'$.
\item\label{item3:dictionnaire} (Voir par exemple le lemme 5
de \cite[\S II.2.1]{Serre83} et l'exemple le suivant, en remarquant
que le faisceau des idéaux au point $P=P_\infty$ (noté $I_P$ dans
loc.~cit.) est $\O_{\bf C}(-P_\infty)$.) Notons $\O_{\bf C}(P_\infty)$
  le fibré en droites sur ${\bf C}$ associé au diviseur
  $P_\infty$. Les conditions suivantes sont équivalentes :

$\bullet$~ le $A$-module projectif $M$ est libre;

$\bullet$~ la restriction de $E^{\overline{M}}$ à ${\bf C}^\circ$ est
triviale;

$\bullet$~ la restriction de $\det (E^{\overline{M}})$ à ${\bf C}^\circ$
est triviale;

$\bullet$~ il existe $z \in\ZZ$ tel que le fibré vectoriel
$\det(E^{\overline{M}})$ soit isomorphe à $\O_{\bf C}
(P_\infty)^{\otimes z}$ (et alors nous avons
$\deg(E^{\overline{M}})=z\,\fff$).
\item\label{item4:dictionnaire} (Voir par exemple
l'affirmation ii) avant la proposition 4 de \cite[\S II.2.1]{Serre83}.)
  Pour tout $z \in\ZZ$, les fibrés $E^{(M,\|\;\|)}\otimes\O_{\bf C}
  (P_\infty)^{\otimes z}$ et $E^{(M,\,q^{-z}\|\;\|)}$ sont isomorphes.
\item\label{item5:dictionnaire} Notons $\overline{M}\oplus
  \overline{M'}$ le $A$-module somme directe $M\oplus M'$ de $M$ et
  $M'$ (donc projectif) muni de la norme du maximum $\|(x,x')\| =
  \max\{\|x\|,\, \|x'\|'\}$ sur le $\wh K$-espace vectoriel $(M\oplus
  M')_{\wh K}={M}_{\wh K} \oplus{M'}_{\wh K}$. Alors
\[
E^{\overline{M}\oplus\overline{M'}}\,=
E^{\overline{M}} \oplus E^{\overline{M'}}\,.
\]
\item\label{item6:dictionnaire} Supposons qu'il existe une $A$-base
  $(b_1,\dots,b_n)$ de $M$ orthogonale pour $\|\;\|$. Notons
  $\overline{M_i}= (A\,b_i, \|\;\|_{\mid \wh K\,b_i})$ le $A$-module
  libre $A\,b_i$ de rang $1$ muni de la restriction de la norme
  $\|\;\|$ de $M_{\wh K}$ à la droite vectorielle $\wh K\,b_i$. Par
  les points \eqref{item3:dictionnaire} et \eqref{item5:dictionnaire},
  nous avons une décomposition $E^{\overline{M}}=\bigoplus_{1\leq i
    \leq n} E^{\overline{M_i}}$ de $E^{\overline{M}}$ en somme directe
  de fibrés en droites sur ${\bf C}$ de restrictions à ${\bf C}^\circ$
  triviales.
\item\label{item7:dictionnaire} (Voir la proposition 4 de
\cite[\S II.2.1]{Serre83}.) Nous dirons que deux fibrés vectoriels sur
  $\bf C$ sont {\it $P_\infty$-stablement isomorphes} s'il existe un
  entier $z\in\ZZ$ tel que les fibrés vectoriels $E\otimes\O_{\bf C}
  (P_\infty)^{\otimes z}$ et $E'$ sur $\bf C$ soient isomorphes. Par
  les points \eqref{item1:dictionnaire}, \eqref{item3:dictionnaire} et
  \eqref{item4:dictionnaire}, l'ensemble des doubles classes
  $\GL_n(\wh\OOO)\bs\GL_n(\wh K)/\GL_n(A)$ s'identifie à l'ensemble
  des classes d'isomorphisme $P_\infty$-stable de fibrés vectoriels de
  rang $n$ sur ${\bf C}$ dont la restriction à ${\bf C}^\circ$ est
  triviale, par l'application qui à la double classe de $g\in\GL_n(\wh
  K)$ associe la classe d'isomorphisme $P_\infty$-stable de
  $E^{(g\,A^n,\;g\|\;\|_{\wh\OOO^n})}$.
\end{enumerate}

\blemm \label{lem:calcdegcovol} Soit $\overline{M}=(M,\|\;\|)$ un
$A$-réseau normé. Nous avons
\[
\deg (E^{\overline{M}})=
-\log_{q_0}\big(\;\overline{\covol} \,(\,\overline{M}\,)\big)\,.
\]
\elemm

Par conséquent, nous avons $\deg (E^{\overline{M}})=0$ si et seulement
si $\overline{M}$ est unimodulaire.
Le membre de droite de cette égalité est l'analogue pour les corps de
fonctions du degré d'Arakelov $-\ln \covol\,(\,\overline{R}\,)$ d'un réseau
euclidien $\overline{R}$, voir \cite[\S 1.3.2]{Bost20}.

\medskip
\dem Soit $(b_1,\dots,b_n)$ une $K$-base de $M_K=M\otimes_AK$
constituée d'éléments de $M$. Notons $M'=\sum_{i=1}^nA\,b_i$, qui est
un $A$-module libre de $A$-base $(b_1,\dots, b_n)$. Remarquons que
$(M')_{K}=M_K$, donc $(M')_{\wh K}=M_{\wh K}$ et $\|\;\|$ est aussi une
norme entière sur $(M')_{\wh K}$.  Notons $\overline{M'}=(M', \|\;\|)$
le $A$-sous-module normé de $\overline{M}$ dont le $A$-module
sous-jacent est $M'$.  Notons pour simplifier $E'= E^{\overline{M'}}$
le sous-fibré vectoriel de $E= E^{\overline{M}}$ associé à
$\overline{M'}$. Soit $(s_1,\dots,s_n)$ une $\O_{{\bf
    C},P_\infty}$-base de la fibre $E'_{P_\infty}$ (qui est un
$\O_{{\bf C},P_\infty}$-module libre). Cette base est orthonormée
pour $\|\;\|$ par la construction de la norme dans l'équivalence de
catégorie $E' \mapsto \overline{M}^{E'}$.  Notons $B$ la matrice de
passage de $(s_1,\dots,s_n)$ à $(b_1,\dots, b_n)$. Nous avons en
particulier $\|\,b_1\wedge\dots \wedge b_n\,\|=|\det B|$.  Puisque
nous avons $M'=\Ga({\bf C}^\circ, E')$, l'élément
$b_1\wedge\dots\wedge b_n=(\det B)\,s_1\wedge\dots\wedge s_n$ est une
section rationnelle du fibré vectoriel $\det (E')$ qui le trivialise
sur ${\bf C}^\circ$ par le point \eqref{item3:dictionnaire} du
dictionnaire. Son diviseur $D$ est donc de support réduit à $P_\infty$
et de plus $D=v(\det B)\,P_\infty$.
Donc
\[
\deg (E')=v(\det B)\;[k_0(P_\infty):k_0]=-\log_{q}(|\det B|)\;\fff
=-\log_{q_0}\|\,b_1\wedge\dots\wedge b_n\,\|\,.
\]
Puisque les normes de $(M')_{\wh K}$ et $M_{\wh K}$ coïncident, par la
construction dans le dictionnaire des fibres en $P_\infty$ de $E'=
E^{\overline{M'}}$ et $E=E^{\overline{M}}$, nous avons
$(E')_{P_\infty}=E_{P_\infty}$. Puisque $E$ et $E'$ ont le même rang
$n$, par le lemme \ref{lem:macri} (et plus précisément par la formule
\eqref{eq:macri}), et puisque les idéaux premiers de l'anneau de
Dedekind $A$ sont les points fermés de ${\bf C}$ différents de
$P_\infty$, nous avons
\begin{align*}
  \deg E-\deg (E')&= \sum_{P\in\P}\dim_{k_0}(E_P/E'_P)=
  \sum_{P\in\P\smallsetminus\{P_\infty\}}\dim_{k_0}(E_P/E'_P)\\&=
  \sum_{P\in\P\smallsetminus\{P_\infty\}}\dim_{k_0}(M_P/M'_P)=
  \dim_{k_0}(M/M')=\log_{q_0} \card(M/M')\,.
\end{align*}
Par la formule \eqref{eq:covolnormwedge}, nous avons par conséquent
\begin{align*}
  \deg E&=\deg (E')+(\deg E-\deg (E'))\\&=
  -\log_{q_0}\|\,b_1\wedge\dots\wedge b_n\,\|
+\log_{q_0}\card(M/M')\\&= -\log_{q_0}\overline{\covol}
\,(\,\overline{M'}\,) +\log_{q_0}\card(M/M')\,.
\end{align*}
Par l'additivé des volumes, nous avons $\frac{\overline{\covol}
  \,(\,\overline{M'}\,)}{\overline{\covol}
  \,(\,\overline{M}\,)}=[M:M']$. Le lemme
\ref{lem:calcdegcovol} en découle. \cqfd

\section{Les pentes de Harder-Narashiman des fibrés
  vectoriels}
\label{sect:pentesHN}

Le contenu de cette partie est essentiellement extrait de
\cite{HarNar75}, auquel nous renvoyons pour les démonstrations. Voir
aussi \cite{Grayson84}, qui prend la convention des signes opposés des
pentes, justifiée par le lemme \ref{lem:calcdegcovol}. Voir enfin
\cite[Appendice A]{Bost20} pour un cadre formel général de
construction des filtrations canoniques de Harder-Narashiman, ainsi
que \cite{Andre09,Chen10b}.

Soit $E$ un fibré vectoriel sur ${\bf C}$ de rang $\rk E=n$. Nous
appellerons {\it affixe} de $E$ le couple $(\rk E,\deg
E)\in\NN\times\ZZ$, et si $\rk E\neq 0$, nous appellerons {\it pente}
de $E$ la pente de la droite vectorielle passant par l'affixe de $E$,
c'est-à-dire le rapport
\[
\mu(E)=\frac{\deg E}{\rk E}\,.
\]
Si $E'$ est un fibré vectoriel sur ${\bf C}$, par la formule
\eqref{eq:degprodtens}, nous avons
\begin{equation}\label{eq:muprodtens}
\mu(E\otimes E')= \mu(E)+\mu(E')\,.
\end{equation}
Si $0\ra A\ra E\ra B\ra 0$ est une suite exacte de fibrés vectoriels
sur ${\bf C}$, par l'additivité des degrés $\deg E=\deg A+\deg B$ et
des rangs $\rk E=\rk A+\rk B$, l'origine $0$ de $\RR^2$ et les affixes
de $A,E,B$ sont les sommets d'un parallélogramme, voir le dessin de
gauche ci-dessous.

\begin{center}
\begin{picture}(0,0)%
\includegraphics{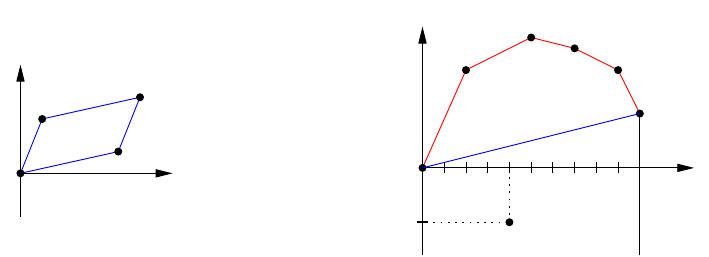}%
\end{picture}%
\setlength{\unitlength}{3812sp}%
\begin{picture}(5855,2114)(11,-1558)
\put(3736,-476){\makebox(0,0)[lb]{\smash{\fontsize{11}{13.2}\usefont{T1}{ptm}{m}{n}{\color[rgb]{1,0,0}$\mu_{\max}(E)$}%
}}}
\put(5221,-136){\makebox(0,0)[lb]{\smash{\fontsize{11}{13.2}\usefont{T1}{ptm}{m}{n}{\color[rgb]{1,0,0}$\mu_{\min}(E)$}%
}}}
\put(4761,-641){\makebox(0,0)[lb]{\smash{\fontsize{11}{13.2}\usefont{T1}{ptm}{m}{n}{\color[rgb]{0,0,1}$\mu(E)$}%
}}}
\put(4366,344){\makebox(0,0)[lb]{\smash{\fontsize{11}{13.2}\usefont{T1}{ptm}{m}{n}{\color[rgb]{0,0,0}$E_2$}%
}}}
\put(4816,209){\makebox(0,0)[lb]{\smash{\fontsize{11}{13.2}\usefont{T1}{ptm}{m}{n}{\color[rgb]{0,0,0}$E_3$}%
}}}
\put(3826,119){\makebox(0,0)[lb]{\smash{\fontsize{11}{13.2}\usefont{T1}{ptm}{m}{n}{\color[rgb]{0,0,0}$E_1$}%
}}}
\put(5356,-421){\makebox(0,0)[lb]{\smash{\fontsize{11}{13.2}\usefont{T1}{ptm}{m}{n}{\color[rgb]{0,0,0}$E_k=E$}%
}}}
\put(3346,-916){\makebox(0,0)[lb]{\smash{\fontsize{11}{13.2}\usefont{T1}{ptm}{m}{n}{\color[rgb]{0,0,0}$0$}%
}}}
\put(1216,-241){\makebox(0,0)[lb]{\smash{\fontsize{11}{13.2}\usefont{T1}{ptm}{m}{n}{\color[rgb]{0,0,0}$E$}%
}}}
\put(1036,-781){\makebox(0,0)[lb]{\smash{\fontsize{11}{13.2}\usefont{T1}{ptm}{m}{n}{\color[rgb]{0,0,0}$B$}%
}}}
\put(276,-356){\makebox(0,0)[lb]{\smash{\fontsize{11}{13.2}\usefont{T1}{ptm}{m}{n}{\color[rgb]{0,0,0}$A$}%
}}}
\put( 26,-921){\makebox(0,0)[lb]{\smash{\fontsize{11}{13.2}\usefont{T1}{ptm}{m}{n}{\color[rgb]{0,0,0}$0$}%
}}}
\put(4276,-1366){\makebox(0,0)[lb]{\smash{\fontsize{11}{13.2}\usefont{T1}{ptm}{m}{n}{\color[rgb]{0,0,0}$F$}%
}}}
\put(5356,-961){\makebox(0,0)[lb]{\smash{\fontsize{11}{13.2}\usefont{T1}{ptm}{m}{n}{\color[rgb]{0,0,0}$n$}%
}}}
\put(3646,-1006){\makebox(0,0)[lb]{\smash{\fontsize{11}{13.2}\usefont{T1}{ptm}{m}{n}{\color[rgb]{0,0,0}$1$}%
}}}
\put(4231,-1006){\makebox(0,0)[lb]{\smash{\fontsize{11}{13.2}\usefont{T1}{ptm}{m}{n}{\color[rgb]{0,0,0}$\rk F$}%
}}}
\put(3016,-1361){\makebox(0,0)[lb]{\smash{\fontsize{11}{13.2}\usefont{T1}{ptm}{m}{n}{\color[rgb]{0,0,0}$\deg F$}%
}}}
\put(5851,-846){\makebox(0,0)[lb]{\smash{\fontsize{11}{13.2}\usefont{T1}{ptm}{m}{n}{\color[rgb]{0,0,0}$\rk$}%
}}}
\put(3376,389){\makebox(0,0)[lb]{\smash{\fontsize{11}{13.2}\usefont{T1}{ptm}{m}{n}{\color[rgb]{0,0,0}$\deg$}%
}}}
\put(1481,-896){\makebox(0,0)[lb]{\smash{\fontsize{11}{13.2}\usefont{T1}{ptm}{m}{n}{\color[rgb]{0,0,0}$\rk$}%
}}}
\put( 81, 64){\makebox(0,0)[lb]{\smash{\fontsize{11}{13.2}\usefont{T1}{ptm}{m}{n}{\color[rgb]{0,0,0}$\deg$}%
}}}
\end{picture}%

\end{center}

Supposons que $\rk E\neq 0$. Le {\it polygone de Harder-Narsimhan} de
$E$ (voir le dessin de droite ci-dessus) est l'enveloppe convexe
fermée $P(E)$ dans $\RR^2$ des affixes de tous les sous-fibrés
vectoriels $F$ de $E$. Par la formule \eqref{eq:satcrois}, prendre
tous les sous-fibrés vectoriels ou seulement les sous-fibrés
vectoriels saturés ne change pas le polygone de Harder-Narsimhan. Il
est contenu dans la bande verticale $[0,n]\times\RR$ et contient
l'ensemble $P^-(E)=\{(x,y)\in[0,n]\times\RR:y\leq \mu(E)\,x\}$ des
points de cette bande en dessous du segment entre l'origine $0$ et
l'affixe de $E$. Les ordonnées des points de $P(E)$ sont
majorées (par exemple car pour tout sous-fibré vectoriel $F$
de $E$ de rang $i\in\llbracket1, n\rrbracket$, par l'inégalité de
Riemann, nous avons $\deg F\leq i(\ggg-1)+h^0({\bf C},F)\leq
n(\ggg-1)+h^0({\bf C},E)$~). L'ensemble des points de $P(E)$
d'ordonnée maximale est une courbe polygonale entre l'origine $0$ et
l'affixe de $E$, graphe d'une fonction concave
$f_E:[0,n]\ra\RR$. Notons que pour tout $i\in\llbracket0,n\rrbracket$,
la valeur $f_E(i)$ est supérieure ou égale (et peut être strictement
supérieure) à la borne supérieure des degrés des sous-fibrés
vectoriels de rang $i$ de $E$. Les pentes successives
\begin{equation}\label{eq:pentesucc}
\mu_i(E)= f_E(i)-f_E(i-1)
\end{equation}
pour $i\in\llbracket1,n\rrbracket$ sont décroissantes
par la convexité de $P(E)$~:
\[
\mu_{\max}(E)=\mu_1(E)\geq\mu_2(E)\geq\dots \geq\mu_n(E)=\mu_{\min}(E)\,.
\]

Notons une différence importante entre les minima successifs
$\lambda_i$ et les pentes $\mu_i$: les premiers sont croissants, les
second décroissants, ce qui justifiera les signes les reliant. De
plus, alors que le premier minimum (la systole) $\lambda_1$ est
atteinte par un sous-$A$-module normé de rang $1$ (même libre), la
première pente $\mu_1$ n'est pas forcément atteinte par un sous-fibré
en droites.

Les {\it pentes de Harder-Narsimhan} de $E$ sont les valeurs
distinctes des pentes successives $\mu_i(E)$ pour
$i\in\llbracket1,n\rrbracket$. Les sommets de $P(E)$ sont les affixes
d'une unique filtration (appelée la {\it filtration de
  Harder-Narsimhan} de $E$) de sous-fibrés vectoriels saturés
\[
0=E_0^{H\!N}\subset E_1^{H\!N}\subset E_2^{H\!N}\subset\dots
\subset E_{k-1}^{H\!N}\subset E_k^{H\!N}=E
\]
où $k\leq n$, et le fibré vectoriel $E$ est dit {\it semi-stable} si
$k=1$, c'est-à-dire si $P(E)$ est réduit à $P^-(E)$. Nous renvoyons
par exemple à \cite[\S 1.3]{HuyLeh10} pour des compléments.

Par exemple, si $E=\oplus_{i=1}^n L_i$ est une somme directe de fibrés
en droites sur ${\bf C}$ de degrés ordonnés de sorte que $\deg
(L_1)\geq\deg (L_2)\geq \dots\geq \deg (L_n)$, alors
\begin{equation}\label{eq:pentesomfibdroit}
\forall\;i\in\llbracket1,n\rrbracket,\qquad \mu_i(E)=\deg (L_i)\,.
\end{equation}

\begin{center}
\begin{picture}(0,0)%
\includegraphics{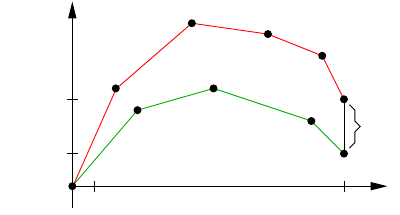}%
\end{picture}%
\setlength{\unitlength}{3812sp}%
\begingroup\makeatletter\ifx\SetFigFont\undefined%
\gdef\SetFigFont#1#2#3#4#5{%
  \reset@font\fontsize{#1}{#2pt}%
  \fontfamily{#3}\fontseries{#4}\fontshape{#5}%
  \selectfont}%
\fi\endgroup%
\begin{picture}(3310,1777)(4711,-1016)
\put(5131,-871){\makebox(0,0)[lb]{\smash{{\SetFigFont{11}{13.2}{\rmdefault}{\mddefault}{\updefault}{\color[rgb]{0,0,0}$0$}%
}}}}
\put(7741,-331){\makebox(0,0)[lb]{\smash{{\SetFigFont{11}{13.2}{\rmdefault}{\mddefault}{\updefault}{\color[rgb]{0,0,0}$\deg E-\deg E'$}%
}}}}
\put(7651,-61){\makebox(0,0)[lb]{\smash{{\SetFigFont{11}{13.2}{\rmdefault}{\mddefault}{\updefault}{\color[rgb]{0,0,0}$E$}%
}}}}
\put(7651,-601){\makebox(0,0)[lb]{\smash{{\SetFigFont{11}{13.2}{\rmdefault}{\mddefault}{\updefault}{\color[rgb]{0,0,0}$E'$}%
}}}}
\put(6301,344){\makebox(0,0)[lb]{\smash{{\SetFigFont{11}{13.2}{\rmdefault}{\mddefault}{\updefault}{\color[rgb]{0,0,0}$P(E)$}%
}}}}
\put(6346,-241){\makebox(0,0)[lb]{\smash{{\SetFigFont{11}{13.2}{\rmdefault}{\mddefault}{\updefault}{\color[rgb]{0,0,0}$P(E')$}%
}}}}
\put(5401,614){\makebox(0,0)[lb]{\smash{{\SetFigFont{11}{13.2}{\rmdefault}{\mddefault}{\updefault}{\color[rgb]{0,0,0}$\deg$}%
}}}}
\put(5446,-961){\makebox(0,0)[lb]{\smash{{\SetFigFont{11}{13.2}{\rmdefault}{\mddefault}{\updefault}{\color[rgb]{0,0,0}$1$}%
}}}}
\put(7606,-916){\makebox(0,0)[lb]{\smash{{\SetFigFont{11}{13.2}{\rmdefault}{\mddefault}{\updefault}{\color[rgb]{0,0,0}$n$}%
}}}}
\put(8006,-816){\makebox(0,0)[lb]{\smash{{\SetFigFont{11}{13.2}{\rmdefault}{\mddefault}{\updefault}{\color[rgb]{0,0,0}$\rk$}%
}}}}
\put(4726,-541){\makebox(0,0)[lb]{\smash{{\SetFigFont{11}{13.2}{\rmdefault}{\mddefault}{\updefault}{\color[rgb]{0,0,0}$\deg E'$}%
}}}}
\put(4771,-106){\makebox(0,0)[lb]{\smash{{\SetFigFont{11}{13.2}{\rmdefault}{\mddefault}{\updefault}{\color[rgb]{0,0,0}$\deg E$}%
}}}}
\end{picture}%

\end{center}

Si $E'$ est un sous-fibré vectoriel de $E$ de même rang $n$, alors
nous avons $P(E')\subset P(E)$. En effet, tout sous-fibré vectoriel
$F'$ de $E'$ de rang $i$ est un sous-faisceau cohérent de $E$, dont le
saturé dans $E$ est un sous-fibré vectoriel $F$ de $E$ de même rang
$i$ tel que $\deg (F')\leq \deg F$. Donc en passant à l'enveloppe
convexe, nous avons $f_{E'}\leq f_E$. En particulier, puisque
$f_E(i)=\sum_{j=1}^i\mu_j(E)$ par la formule \eqref{eq:pentesucc},
nous avons $\mu_1(E')\leq \mu_1(E)$,
\[
\forall\;i\in\llbracket 1,n\rrbracket,\quad
\sum_{j=1}^i\mu_j(E')\leq \sum_{j=1}^i\mu_j(E)\,,
\]
et
\[
\sum_{j=1}^n\mu_j(E')=f_{E'}(n)= f_{E}(n) -( \deg E-\deg (E'))=
\sum_{j=1}^n\mu_j(E) -(\deg E-\deg (E'))\;.
\]
  
\blemm Pour tout $i\in\llbracket1,n\rrbracket$ et tout fibré en
droites $L$ sur ${\bf C}$, nous avons
\begin{equation}\label{eq:mumaxprodtensdroite}
\mu_i(E\otimes L)= \mu_i(E)+\deg L\,.
\end{equation}
\elemm

\dem Pour tout fibré vectoriel $E'$ de rang $n$ sur ${\bf C}$ et pour
tout $i\in\llbracket1,n\rrbracket$, notons $\operatorname{degmax}_i
(E')$ le degré maximum d'un sous-fibré vectoriel de rang $i$ de $E'$.
Pour tout sous-fibré vectoriel $F$ de $E$ de rang $i$ de degré maximal
(donc saturé), le fibré vectoriel $F\otimes L$ est un sous-fibré
vectoriel de rang $i$ de $E\otimes L$, et $\deg (F\otimes L)=\deg
F+i\deg L$.  Donc nous avons $\operatorname{degmax}_i (E\otimes L)\geq
\operatorname{degmax}_i (E) +i\deg L$. Par la linéarité de
l'application $t\mapsto t\deg L$ et par passage à l'enveloppe convexe,
nous avons donc $f_{E\otimes L}(i) \geq f_{E}(i)+i\deg L$. En
remplaçant dans cette inégalité $(E,L)$ par $(E\otimes L,
L^{\widecheck{\;}}\;)$ où $L^{\widecheck{\;}}$ est le dual de $L$, et
puisque $\deg (L^{\widecheck{\;}})=-\deg L$, nous avons $f_{E\otimes
  L}(i) = f_{E}(i)+i\deg L$. Le résultat découle donc de la formule
\eqref{eq:pentesucc} par différence.  \cqfd

\section{Pentes et minima}
\label{sect:pentmin}

\medskip\noindent{\bf Relation entre pente maximale et systole. }
Soient $n\in\NN\ssm\{0\}$, $\overline{M}= (M,\|\cdot\|)$
un $A$-réseau normé de rang $n$ dont la norme est entière, et $E=
E^{\overline{M}}$ son fibré vectoriel sur ${\bf C}$ associé.

\blemm \label{lem:relatpentmaxsys} Nous avons
$\displaystyle\lambda_1(\,\overline{M}\,)\in
\Big\llbracket\,\Big\lceil\frac{-\mu_{\max}(E)}{\fff} \Big\rceil,
\Big\lfloor\frac{\ggg-\mu_{\max}(E)-1}{\fff} \Big\rfloor+1\,\Big\rrbracket$.
\elemm

En particulier, si $\ggg=0$ et $\fff=1$, alors
$\lambda_1(\,\overline{M}\,)= -\,\mu_{\max}(E^{\overline{M}})$.

\medskip \dem Montrons d'abord la borne supérieure. Soit $F$ un
sous-fibré vectoriel non nul (forcément saturé) de $E$ tel que
$\mu_{\max}(E)=\mu(F)$.  Par le dictionnaire, il existe un
sous-$A$-module normé $\overline{N}= (N, \|\cdot\|_{\mid N_{\wh K}})$
de $\overline{M}$ tel que le sous-fibré vectoriel $F$ soit isomorphe à
$E^{\overline{N}}$.  Posons $z =
\Big\lfloor\frac{\ggg-\mu_{\max}(E)-1}{\fff} \Big\rfloor+1$ et
$F'=F\otimes \O_{\bf C} (P_\infty)^{\otimes z}$. Par la formule
\eqref{eq:muprodtens}, puisque $\O_{\bf C} (P_\infty)^{\otimes z}$ est
un fibré en droites sur $\bf C$ de degré $\fff z$ et puisque $\lfloor
t\rfloor+1>t$ pour tout $t\in\RR$, nous avons
\[
\mu(F') =\mu(F)+\fff \,z> \mu(F)+\ggg-\mu_{\max}(E)-1=\ggg-1\;.
\]
Par l'inégalité de Riemann, nous avons donc
\[
h^0({\bf C},F')\geq \deg(F')+\rk(F')(1-\ggg)=
\rk(F') \big(\mu(F')+1-\ggg\big)>0\;.
\]
Par le point \eqref{item4:dictionnaire} du dictionnaire et la formule
\eqref{eq:relatsecttotal}, nous avons donc
\[
\Ga({\bf C},F') =\{s\in N:q^{-z}\|s\|\leq 1\}\neq \{0\}\;,
\]
donc $\lambda_1(\,\overline{M}\,)\leq \lambda_1(\,\overline{N}\,) \leq
z$.

Montrons maintenant la borne inférieure. Soit $v\in M$ un plus court
vecteur non nul de $\overline{M}$, de sorte que $\|v\| =
q^{\lambda_1(\,\overline{M}\,)}$. Notons $\overline{N} = (N=Av,
\|\;\|_{\mid N_{\wh K}})$, qui est un sous-$A$-réseau normé de
$\overline{M}$ de rang $1$ dont la norme est entière, et
$F=E^{\overline{N}}$ le sous-fibré vectoriel de $E$ de rang $1$
associé. Alors par le lemme \ref{lem:calcdegcovol} et la formule
\eqref{eq:covolnormwedge} appliquée à $\Lambda=A\,v$ et
$n=1$, nous avons
\begin{align*}
\mu_{\max} (E)\geq \mu(F)=\deg F &=-\log_{q_0}\big(\;\overline{\covol}\,
(\,\overline{N}\,)\big)=-\log_{q_0}\|v\|\\ &=-\,\fff\log_q\|v\|=
  -\,\fff\,\lambda_1(\,\overline{M}\,)\;.
\end{align*}
Donc $\lambda_1(\,\overline{M}\,) \geq -\frac{\mu_{\max}(E)}{\fff}$,
ce qui montre le résultat, sachant que
$\lambda_1(\,\overline{M}\,) \in\ZZ$ puisque la norme $\|\;\|$ est entière.
\cqfd

\medskip\noindent{\bf Quasi-décalage des pentes.}  Soient $E$ un fibré
vectoriel sur ${\bf C}$ de rang $n$, et $L$ un sous-fibré en droites
de $E$ de degré maximal, donc saturé.  

\blemm \label{lem:quasidecal}
Pour tout $i\in\llbracket2,n\rrbracket$, nous avons
\[
 |\;\mu_{i-1}(E/L)- \mu_{i}(E)\;|\leq 2\,\ggg+2\fff-2\,.
\]
\elemm

Pour $i=2$, la démonstration (voir les formules
\eqref{eq:majoEmodLparE} et \eqref{eq:majoEparEmodL} pour $i=2$)
montre que nous avons la légère amélioration
\begin{equation}\label{eq:decalageunmodL}
|\;\mu_{1}(E/L)- \mu_{2}(E)\;|\leq\ggg+\fff-1\,.
\end{equation}

\medskip
\dem Montrons tout d'abord que nous avons
\begin{equation}\label{eq:droitdegmax}
\mu_{\max}(E)-\ggg-\fff+1 \leq \deg L=\mu(L)\leq\mu_{\max}(E)\,.
\end{equation}
Notons $\overline{M}^{\,E}=(M,\|\;\|)$ le $A$-réseau normé associé à
$E$. Soit $b$ un élément de $M$ tel que $\|b\|=
q^{\lambda_1(\,\overline{M}^{\,E}\,)}$. Soit $\overline{M'}=(A\,b,
\|\;\|_{\mid \wh K\,b})$ le sous-$A$-réseau normé de
$\overline{M}^{\,E}$ engendré par $b$, et $E^{\overline{M'}}$ le
sous-fibré en droites de $E$ associé.
Par la formule \eqref{eq:covolnormwedge} appliquée à $\Lambda=A\,b$ et
$n=1$, nous avons $\overline{\covol}\,(\,\overline{M'}\,)=\|b\|$.  Par
les lemmes \ref{lem:relatpentmaxsys} et \ref{lem:calcdegcovol}, nous
avons alors
\begin{align*}
  \mu_{\max}(E)\,-\,\ggg -\fff+1&=
  -\fff\Big(\frac{\ggg-\mu_{\max}(E)-1}{\fff}+1\Big) \leq
  -\fff\Big(\Big\lfloor\frac{\ggg-\mu_{\max}(E)-1}{\fff}\Big\rfloor+1\Big)
  \\&\leq -\,\fff\,\lambda_1(\,\overline{M}^{\,E}\,)
=-\,\fff\,\log_q\|b\|=
-\log_{q_0}\|b\|\\&=-\log_{q_0}\overline{\covol}\,(\,\overline{M'}\,)=
\deg (E^{\overline{M'}})\,.
\end{align*}
Ceci démontre donc la formule \eqref{eq:droitdegmax} par la maximalité
du degré de $L$.

Maintenant, soit $i\in\llbracket2,n\rrbracket$. Notons $p: E\ra E/L$
la projection canonique.  Pour tout sous-fibré vectoriel $F$ de $E/L$
de rang $i-1$, l'image réciproque $p^{-1}(F)$ est un sous-fibré
vectoriel de $E$ de rang $i$ et $L$ est saturé dans $p^{-1}(F)$.  Par
l'additivité des degrés dans la suite exacte de fibrés vectoriels
$0\ra L\ra p^{-1}(F)\ra F\ra 0$, nous avons $\deg F=\deg
(p^{-1}(F))-\deg L$.  En passant à l'enveloppe convexe et en rappelant
que $f_E$ désigne la fonction concave dont le graphe est le bord
supérieur du polygone de Harder-Narasimhan de $E$, nous avons donc
\[
f_{E/L}(i-1)\leq f_E(i)-\deg L\,.
\]
Puisque $f_E(i)=\sum_{j=1}^i\mu_j(E)$ par la formule
\eqref{eq:pentesucc}, et puisque 
\begin{equation}\label{eq:minodegvectdroitmax}
\deg L\geq \mu_1(E)-\ggg-\fff+1
\end{equation}
par la
formule \eqref{eq:droitdegmax}, nous avons donc
\begin{equation}\label{eq:majoEmodLparE}
\sum_{j=2}^i\mu_{j-1}(E/L)\leq\sum_{j=2}^i\mu_{j}(E)+\ggg+\fff-1\,.
\end{equation}

Réciproquement, soient $F$ un sous-fibré vectoriel de $E$ de rang $i$,
et $F'$ le sous-fibré vectoriel de $E/L$ saturé du sous-faisceau
$p(F)$. Si $F$ contient $L$, alors $L$ est saturé dans $F$, nous avons
$p(F)=F/L$, le fibré vectoriel $F'$ est de rang $i-1$ et
\begin{align*}
  \deg F&=\deg L+\deg (F/L)\leq \deg L+\deg (F') \\&\leq \mu_1(E)+
  f_{E/L}(i-1)=\mu_{1}(E)+\sum_{j=2}^i\mu_{j-1}(E/L)\,.
\end{align*}
Si $F$ ne contient pas $L$, alors $F'$ est de rang $i$. Par la
formule \eqref{eq:majoEmodLparE} avec $i=2$ et par la décroissance des
pentes de Harder-Narashiman, nous avons
\begin{align*}
  \deg F&\leq \deg (F')\leq f_{E/L}(i)= \mu_{1}(E/L)+ \sum_{j=2}^i\mu_{j}(E/L)
  \\&\leq (\mu_{2}(E)+\ggg+\fff-1)+ \sum_{j=2}^i\mu_{j}(E/L)
  \leq \mu_{1}(E)+\ggg+\fff-1+ \sum_{j=2}^i\mu_{j-1}(E/L)\,.
\end{align*}
Par passage à l'enveloppe convexe, nous avons donc
\[
f_E(i)\leq\max\{\mu_1(E)+ f_{E/L}(i-1),\;f_{E/L}(i)\}
\leq\mu_{1}(E)+\ggg+\fff-1+ \sum_{j=2}^i\mu_{j-1}(E/L)\,.
\]
En enlevant $\mu_{1}(E)$ des deux côtés de cette inégalité, nous avons
\begin{equation}\label{eq:majoEparEmodL}
\sum_{j=2}^i\mu_{j}(E)\leq \sum_{j=2}^i\mu_{j-1}(E/L)+\ggg+\fff-1\,.
\end{equation}
Avec la convention usuelle des sommes vides, les formules
\eqref{eq:majoEparEmodL} et \eqref{eq:majoEmodLparE} donnent
\begin{align*}
\mu_i(E)&=\sum_{j=2}^i\mu_{j}(E)-\sum_{j=2}^{i-1}\mu_{j}(E)\\&\leq
\big(\sum_{j=2}^i\mu_{j-1}(E/L)+\ggg+\fff-1\big) -
\big(\sum_{j=2}^{i-1}\mu_{j-1}(E/L)-\ggg-\fff+1\big)
\\&=\mu_{i-1}(E/L)+2\ggg+2\fff-2\,,
\end{align*}
et de même $\mu_{i-1}(E/L)\leq\mu_i(E)+2\ggg+2\fff-2$.  Le résultat en
découle.
\cqfd

\medskip\noindent{\bf Rappels sur les modules de type fini sur les
  anneaux de valuation discrète.} Soient $\A$ un anneau de valuation
discrète, $v$ sa valuation, $\M$ son idéal maximal, $\varpi$ une
uniformisante de sorte que $\M=\varpi\A$, $|\cdot|=q^{-v(\cdot)}$ sa
valeur absolue normalisée par $|\varpi|=q^{-1}$, $\K$ le corps des
fractions de $\A$, et $T,T',T''$ trois $\A$-modules de type fini.

Supposons que les $\A$-modules $T$ et $T'$ soient sans torsion. Ils
sont donc libres. Nous munissons le $\K$-espace vectoriel
$T_\K=T\otimes_A\K$ (respectivement $T'_\K=T'\otimes_A\K$) de la norme
du maximum dans n'importe quelle $\A$-base de $T$ (respectivement
$T'$), et le $\K$-espace vectoriel $\Hom_\K(T_\K,T'_\K)$ de la norme
d'opérateur $\interleave\,\cdot\,\interleave$ correspondante.

Supposons que le $\A$-module $T''$ soit de torsion, donc isomorphe à
$\prod_{i=1}^n\A/\M^{a_i}$ où $n\in\NN$ et $a_1,\dots, a_n\in\NN\ssm
\{0\}$ (uniques si ordonnés décroissants). Rappelons que la longueur
de $T''$ est $\operatorname{lg} \,(T'')=\sum_{i=1}^n a_i$ et par
définition, l'{\it ordre d'annulation} $\operatorname{ann}\,(T'')$ de
$T''$ est la borne inférieure des entiers $a\in\NN$ tels que
$\varpi^{a}T''=\{0\}$, de sorte que $\operatorname{ann}\,(T'')
=\max_{1\leq i\leq n} a_i$. Notons qu'alors
\[
\operatorname{lg}\,(T'')\geq \operatorname{ann}\,(T'')\;.
\]

\blemm \label{lem:controlinverselong} Soit $\psi:T\ra T'$ un morphisme
de $\A$-modules de type fini. Supposons $T$ et $T'$ sans torsion.
L'extension $\K$-linéaire $\psi_\K:T_\K\ra T'_\K$ de $\psi$ est un
isomorphisme de $\K$-espaces vectoriels si et seulement si $\psi$ est
injectif et si le $\A$-module $T'/\operatorname{im}(\psi)$ est de
torsion, et alors nous avons
$\interleave\,(\psi_\K)^{-1}\,\interleave=
q^{\operatorname{ann}\,(T'/\operatorname{im}(\psi))}\leq
q^{\operatorname{lg}\,(T'/\operatorname{im}(\psi))}$.  \elemm

\dem L'équivalence étant immédiate, montrons la dernière affirmation.
Le morphisme d'extension $\psi\mapsto \psi_\K$
de $\Hom_\A(T,T')$ dans $\Hom_\K(T_\K,T'_\K)$ est injectif. Son image
est égale à $\{\psi'\in\Hom_\K(T_\K,T'_\K): \interleave\,\psi'\,
\interleave\leq 1\}$ car les boules unités des normes de $T_\K$ et
$T'_\K$ sont exactement $T$ et $T'$.

Notons $\mu:T'\ra T'$ le morphisme de $\A$-modules de multiplication
par $\varpi^{\operatorname{ann}\,(T'/\operatorname{im}(\psi))}$. Son
image est contenue dans l'image $\operatorname{im}(\psi)$ de $\psi$,
par la définition de l'ordre d'annulation. Donc il peut être
postcomposé par le morphisme inverse de l'isomorphisme de $\A$-modules
$\psi:T\ra\psi(T)$.  Puisque la composition $\psi^{-1}\circ\mu$ est un
morphisme de $\A$-modules, son extension $\K$-linéai\-re, qui est
$\varpi^{\operatorname{ann}\, (T'/\operatorname{im}(\psi))}
(\psi_\K)^{-1}$, est de norme d'opérateur au plus $1$ par ce qui
précède. Ceci montre que nous avons $\interleave\,(\psi_\K)^{-1}\,
\interleave\leq q^{\operatorname{ann} (T'/\operatorname{im}(\psi))}$,
et l'inégalité réciproque (dont nous n'aurons pas besoin) est laissée
au lecteur.  \cqfd

\medskip\noindent{\bf Quasi-scindage de fibrés vectoriels sur ${\bf
    C}$.}  Lorsque ${\bf C}$ est la droite projective, il est bien
connu (voir par exemple \cite{Grothendieck57} dont la démonstration
est en fait indépendante du corps de base, ainsi que
\cite[Theo.~1.3.1]{HuyLeh10}) que tout fibré vectoriel sur ${\bf C}$
est somme directe de fibrés en droites. Le résultat suivant donne une
estimation optimale du défaut d'un tel scindage.

\btheo \label{theo:quasiscind} Pour tout fibré vectoriel $E$ sur ${\bf
  C}$ de rang $n$, il existe des sous-fibrés vectoriels $L_1,\dots,
L_n$ de $E$ de rang $1$ tels que le morphisme de faisceau somme
$\Sigma:L_1\oplus \dots\oplus L_n\ra E$ soit un isomorphisme sur ${\bf
  C}^\circ$ et que le $\O_{{\bf C},P_\infty}$-module de type fini de
torsion $E_{P_\infty}/\operatorname{im}(\Sigma)_{P_\infty}$ vérifie
\begin{equation}\label{eq:quasiscind}
\fff\,\operatorname{lg}(E_{P_\infty}/\operatorname{im}(\Sigma)_{P_\infty})=
\deg E-\sum_{i=1}^n\deg (L_i)\leq \frac{n(n-1)}{2}\,(4\ggg+3\fff-3)\,.
\end{equation}
\etheo

Lorsque $\ggg=0$ et $\fff=1$, nous retrouvons le fait que $E$ est
somme directe de sous-fibrés en droites, voir
\cite[Theo.~1.3.1]{HuyLeh10}).

\medskip
\dem L'égalité à droite dans la formule \eqref{eq:quasiscind}, qui
utilise le fait que $\Sigma$ soit un isomorphisme sur ${\bf C}^\circ$,
découle du lemme \ref{lem:macri}.  La démonstration de l'inégalité à
gauche dans la formule \eqref{eq:quasiscind} procède par récurrence
sur $n$. Le résultat si $n=1$ est immédiat, avec $L_1=E$, $\Sigma=\id$
et donc $E=\operatorname{im}(\Sigma)$.  L'outil clef est le lemme
suivant.

\blemm \label{lem:quasiscindk} Soient $S,E,Q$ des fibrés vectoriels
sur ${\bf C}$ et $\E$ une suite exacte de faisceaux $0\ra
S\stackrel{\iota}{\longrightarrow} E\stackrel{p}{\longrightarrow} Q\ra
0$. Si
\[
k=\max\Big\{0,\Big\lfloor\frac{1}{\fff}\,(2\ggg-2-\deg S+\mu_{\max}(Q))
\Big\rfloor+1\Big\}
\]
et si $S$ est de rang $1$, alors il existe un scindage
$s^\circ:Q_{\mid{\bf C}^\circ}\ra E_{\mid{\bf C}^\circ}$ de $\E$
au-dessus de ${\bf C}^\circ$ de pôle d'ordre au plus $k$ en
$P_\infty$, c'est-à-dire telle que la post-composition par l'inclusion
$E\ra E\otimes\O_C (P_\infty)^{\otimes k}$ de l'extension sur ${\bf
  C}$ de $s^\circ$ soit régulière.
\elemm

\dem La suite exacte $0\ra S_{\mid{\bf C}^\circ}
\stackrel{\iota}{\longrightarrow} Q_{\mid{\bf C}^\circ}
\stackrel{p}{\longrightarrow} E_{\mid{\bf C}^\circ}\ra 0$ scinde car
la variété ${\bf C}^\circ$ est affine.

Notons $j:Q(-kP_\infty)=Q\otimes\O_{\bf C}(P_\infty)^{\otimes-k} \ra
Q$ l'inclusion dans $Q$ de son sous-faisceau saturé de ses éléments
ayant un pôle d'ordre au plus $k$ en $P_\infty$. Notons $E_k$
l'extension de $Q(-kP_\infty)$ par $S$ tirée en arrière de $\E$ par
$j$, et $\E_k$ la suite exacte de faisceaux associée, de sorte que
nous avons un diagramme commutatif de ligne exactes
\[\begin{array}{ccccccc}
0 \ra & S & \longrightarrow & E &\stackrel{p}{\longrightarrow}
&Q&\ra 0{\;}\\ & \| & & \uparrow &
&\;\uparrow_j &\\
0 \ra & S & \longrightarrow & E_k &\stackrel{p_k}{\longrightarrow}
&Q(-kP_\infty)&\ra 0\,.
\end{array}
\]
La conclusion du lemme est vérifiée si et seulement s'il existe un
morphisme de faisceaux $s'_k:Q(-kP_\infty)\ra E$ tel que $p\circ
s'_k=j$, donc si et seulement s'il existe un morphisme de faisceaux
$s_k:Q(-kP_\infty)\ra E_k$ tel que $p_k\circ s_k=\id_{Q(-kP_\infty)}$,
donc si et seulement si la classe de $\E_k$ dans
$\operatorname{Ext}^1_{\O_{\bf C}}(Q(-kP_\infty),S)\simeq H^1({\bf C},
Q^{\;\widecheck{\;}}\otimes S \otimes\O_{\bf C} (P_\infty)^{\otimes
  k})$ est nulle.

Puisque le degré du fibré en droites canonique $\Omega^1_{\bf C}$ est
$\deg (\Omega^1_{\bf C})=2\,\ggg-2$, puisque $\deg(\O_{\bf C}
(P_\infty) ^{\otimes -k})=-\,k\,\fff$, par la formule
\eqref{eq:mumaxprodtensdroite} car $S$ est de rang $1$, par la
définition de $k$, et puisque $\lfloor t\rfloor+1>t$ pour tout
$t\in\RR$, nous avons
\begin{align*}
\mu_{\max}\big(Q\otimes S^{\;\widecheck{\;}}
\otimes\O_{\bf C}(P_\infty)^{\otimes -k}\otimes\Omega^1_{\bf C}\big)=
\mu_{\max}(Q)-\deg S -\,k\,\fff+2\,\ggg -2<0\,.
\end{align*}
Puisque sa première pente est strictement négative et que les pentes
sont décroissantes, tous les degrés des sous-fibrés vectoriels non
triviaux de $Q\otimes S^{\;\widecheck{\;}}\otimes \O_{\bf C}
(P_\infty)^{\otimes -k}\otimes\Omega^1_{\bf C}$ sont strictement
négatifs, donc ce fibré n'a pas de section régulière non nulle. Par la
dualité de Serre, nous avons par conséquent
\[
h^1({\bf C},Q^{\;\widecheck{\;}}\otimes S
\otimes\O_{\bf C}(P_\infty)^{\otimes k})=h^0({\bf C},Q\otimes S^{\;\widecheck{\;}}
\otimes\O_{\bf C}(P_\infty)^{\otimes -k}\otimes\Omega^1_{\bf C})=0\,.
\]
Le résultat en découle.
\cqfd

\medskip
Supposons maintenant que $n\geq 2$ et soit $L_1$ un sous-fibré en
droites de $E$ de degré maximal, donc saturé. Par récurrence, il
existe des sous-fibrés en droites $L'_2,\dots, L'_n$ du fibré
vectoriel $E/L_1$ tels que le morphisme somme $\Sigma':L'_2\oplus
\dots\oplus L'_n\ra E/L_1$ soit un isomorphisme sur ${\bf C}^\circ$ et
que
\[
\fff\,\operatorname{lg}((E/L_1)_{P_\infty}/
\operatorname{im}(\Sigma')_{P_\infty})=
\deg (E/L_1)-\sum_{i=2}^n\deg (L'_i)\leq
\frac{(n-1)(n-2)}{2}\,(4\ggg+3\fff-3)\,.
\]
Notons $k\in\NN$ l'entier défini dans le lemme \ref{lem:quasiscindk}
appliqué à l'inclusion $\iota:S=L_1\ra E$ et à la projection canonique
$p:E\ra Q=E/L_1$ (de sorte que $S$ est bien de rang $1$).  Par les
formules \eqref{eq:minodegvectdroitmax} et \eqref{eq:decalageunmodL}
et par la décroissance des pentes, nous avons
\begin{align*}
  k&=\max\Big\{0,\Big\lfloor\frac{1}{\fff}\,
  (2\ggg-2-\deg (L_1)+\mu_{\max}(E/L_1))\Big\rfloor+1\Big\}\\&
  \leq \max\Big\{0,\frac{1}{\fff}\,
  (2\ggg-2-(\mu_1(E)-\ggg-\fff+1)+(\mu_2(E)+\ggg+\fff-1))+1\Big\}
  \\&\leq\max\Big\{0,\frac{1}{\fff}\,(4\,\ggg+3\fff-4)\Big\}\leq
  \frac{1}{\fff}\,(4\,\ggg+3\fff-3)\,.
\end{align*}
Pour tout $i\in\llbracket 2,n\rrbracket$, soit
$L_i$ un sous-fibré en droite de $E$ tel que $(L_i)_{\mid {\bf C}
  ^\circ} =s^0((L'_i)_{\mid {\bf C}^\circ})$.  Par le lemme
\ref{lem:quasiscindk}, l'ordre d'annulation, en tant que $\O_{{\bf C}
  ,P_\infty}$-module de torsion monogène, de la fibre en $P_\infty$ du
faisceau quotient $L'_i/p(L_i)$ est au plus $k$. Donc, par le lemme
\ref{lem:macri}, nous avons
\begin{align*}
\deg (L'_i)-\deg (L_i)&=\fff\,\operatorname{lg}
\big((L'_i)_{P_\infty}/(p(L_i))_{P_\infty}\big)=
\fff\,\operatorname{ann}\big((L'_i)_{P_\infty}/(p(L_i))_{P_\infty}\big)
\\&\leq \fff\, k\leq 4\,\ggg+3\fff-3\,.
\end{align*}
Le morphisme somme $\Sigma:L_1\oplus \dots\oplus L_n\ra E$ est un
isomorphisme sur ${\bf C}^\circ$ et nous avons
\begin{align*}
\fff\,\operatorname{lg}((E_{P_\infty}/\operatorname{im}(\Sigma)_{P_\infty})&=
\deg E-\sum_{i=1}^n\deg (L_i)\\&= \deg (E/L_1)-\sum_{i=2}^n\deg (L'_i) +
\sum_{i=2}^n(\deg (L'_i)-\deg (L_i))\\&
\leq \frac{(n-1)(n-2)}{2}\,(4\ggg+3\fff-3)+(n-1)\,(4\,\ggg+3\fff-3)
\\&=\frac{n(n-1)}{2}\,(4\ggg+3\fff-3)\,.
\end{align*}
Le résultat en découle.
\cqfd

\bcoro \label{coro:pentdegdroit}
Pour tout fibré vectoriel $E$ sur ${\bf C}$ de rang $n$, il existe des
sous-fibrés vectoriels $L_1,\dots, L_n$ de $E$ de rang $1$ tels que le
morphisme de faisceau somme $L_1\oplus \dots\oplus L_n\ra E$ soit un
isomorphisme sur ${\bf C}^\circ$ et pour tout $i\in\llbracket 1,
n\rrbracket$, nous ayons
\[
|\,\deg (L_i)-\mu_i(E)\,|\leq (n-1)\,(6\ggg+5\fff-5)\,.
\]
\ecoro

\dem Posons $c_n= (n-1)\,(6\ggg+5\fff -5)$.  Montrons le résultat par
récurrence sur $n$. C'est immédiat si $n=1$, en prenant $L_1=E$.
Supposons donc $n\geq 2$. Soit $L_1$ un sous-fibré en droites de $E$
de degré maximal, donc saturé. Par la formule
\eqref{eq:minodegvectdroitmax} et puisque $n\geq 2$, nous avons
\[
|\,\deg (L_1)-\mu_1(E)\,|\leq \ggg+\fff-1\leq c_{n}\,.
\]
Par récurrence, il existe des sous-fibrés vectoriels $L'_2,\dots,
L'_n$ de rang $1$ du fibré vectoriel $E/L_1$ tels que le morphisme
somme $L'_2\oplus \dots\oplus L'_n\ra E/L_1$ soit un isomorphisme sur
${\bf C}^\circ$ et pour tout $i\in\llbracket 2, n\rrbracket$, nous
ayons
\[
|\,\deg (L'_i)-\mu_{i-1}(E/L_1)\,|\leq c_{n-1}\,.
\]
Comme dans la démonstration du théorème \ref{theo:quasiscind}, il
existe des sous-fibrés vectoriels $L_2,\dots, L_n$ de $E$ de rang $1$
tels que le morphisme somme $L_1\oplus \dots\oplus L_n\ra E$ soit un
isomorphisme sur ${\bf C}^\circ$ et pour tout $i\in\llbracket 2,
n\rrbracket$, nous ayons
\[
|\,\deg (L'_i)-\deg (L_i)\,|\leq 4\ggg+3\fff-3\,.
\]
Donc pour tout $i\in\llbracket 2, n\rrbracket$, par le lemme
\ref{lem:quasidecal}, nous avons
\begin{align*}
  &|\,\deg (L_i)-\mu_i(E)\,|\\=\;&|\,(\deg (L_i)-\,\deg (L'_i))+
  (\deg (L'_i)-\mu_{i-1}(E/L_1))+(\mu_{i-1}(E/L_1)-\mu_i(E))\,|
  \\\leq\;& (4\,\ggg+3\fff-3) +|\,\deg (L'_i)-\mu_{i-1}(E/L_1)\,|+
  (2\,\ggg+2\fff-2)\\ \leq\;& c_{n-1} + (6\ggg+5\fff-5)  =c_n\;.
\end{align*}
Le résultat en découle.
\cqfd

\medskip\noindent{\bf Relation entre pentes et minima successifs. }
Soient $n\in\NN\ssm\{0\}$, $\overline{M}= (M,\|\cdot\|)$ un $A$-réseau
normé de rang $n$ de norme entière, et $E=E^{\overline{M}}$ son fibré
vectoriel sur ${\bf C}$ associé. Le résultat suivant compare les
minima successifs $\lambda_1 (\,\overline{M}\,)\leq \lambda_2
(\,\overline{M}\,)\leq \dots\leq\lambda_n(\,\overline{M}\,)$ de
$\overline{M}$ aux pentes de Harder-Narashiman $\mu_1 (E)\geq
\mu_2(E)\geq \dots\geq\mu_n(E)$ de $E$.

\btheo\label{theo:relatpenminsucc} Il existe une permutation $\sigma$
de $\llbracket 1,n\rrbracket$ telle que pour tout $i\in\llbracket
1,n\rrbracket$, nous ayons
\[
\textstyle\lambda_i(\,\overline{M}\,)\in \Big\llbracket\,
\Big\lceil\frac{-\,\mu_{\sigma(i)}(E^{\overline{M}})\,-\,
(n-1)(6\ggg+5\fff-5)\,-\,\frac{n(n-1)}{2}(4\ggg+3\fff-3)}{\fff}\Big\rceil,
\Big\lfloor\frac{-\,\mu_{\sigma(i)}(E^{\overline{M}})\,+\,(n-1)(6\ggg+5\fff-5)
  \,+\,\ggg-1}{\fff}\Big\rfloor+1\,\Big\rrbracket\;.
\]
\etheo

En particulier, si $\ggg=0$ et $\fff=1$, alors par les propriétés de
croissance des $\lambda_i$ et de décroissance des $\mu_i$ pour
$i\in\llbracket 1,n\rrbracket$, nous pouvons prendre $\sigma=\id$ de
sorte que nous ayons $\lambda_i(\,\overline{M}\,)=-
\mu_{i}(E^{\overline{M}})$ pour tout $i\in\llbracket 1,n\rrbracket$.
Nous ne savons pas si les constantes qui apparaissent ci-dessus sont
optimales pour $(\ggg,\fff)$ quelconque. En utilisant le fait que les
suites $(\lambda_i(\,\overline{M}\,))_{1\leq i\leq n}$ et $(-\,\mu_{j}
(E^{\overline{M}}))_{1\leq j\leq n}$ sont croissantes, et en posant
\[
\textstyle c'_{\ggg,\fff,n}=\max\big\{\big\lfloor\frac{\,
(n-1)(6\ggg+5\fff-5)\,+\,\frac{n(n-1)}{2}(4\ggg+3\fff-3)}{\fff}\big\rfloor
,\big\lceil\frac{(n-1)(6\ggg+5\fff-5)
  \,+\,\ggg-1}{\fff}\big\rceil+1\big\}\,,
\]
la formule \eqref{eq:relatpenminsuccintro} du théorème
\ref{theo:mainintro} de l'introduction en découle. Remarquons que
$c'_{0,1,n}=0$.

\medskip \dem Par le corollaire \ref{coro:pentdegdroit} appliqué à $E=
E^{\overline{M}}$, soient $L_1,\dots, L_n$ des sous-fibrés vectoriels
de $E$ de rang $1$ tels que le morphisme somme $\varphi: L_1\oplus
\dots\oplus L_n\ra E$ soit un isomorphisme sur ${\bf C}^\circ$ et
$|\,\deg (L_i)-\mu_i(E)\,|\leq (n-1)\,(6\ggg+5\fff-5)$ pour tout $i\in
\llbracket 1, n\rrbracket$. Par l'équivalence de catégorie, soient
$\overline{M}_1=(M_1,\|\;\|_1),\dots, \overline{M}_n= (M_n,\|\;\|_n)$
des $A$-réseaux normés de rang $1$ de norme entière tels que
$E^{\overline{M}_i}=L_i$ pour tout $i\in\llbracket
1,n\rrbracket$.

Puisque $\varphi$ est un isomorphisme sur ${\bf C}^\circ$, nous
pouvons supposer que $M_1,\dots, M_n$ sont des sous-$A$-modules de $M$
tels que le morphisme associé (voir le début de la partie
\ref{sect:fibvect}) $\phi=\phi(\varphi):M_1\oplus\dots\oplus M_n\ra M$
soit un isomorphisme de $A$-modules.  Notons
\[
\overline{M'}=\big(M'=M_1\oplus\dots\oplus M_n,
\|\;\|':(x_1,\dots,x_n)\mapsto \max_{1\leq i\leq n}\|x_i\|_i\big)
\]
le $A$-module normé somme directe orthogonale. Soit $\sigma$ une
permutation de $\llbracket 1,n\rrbracket$ telle que
\begin{equation}\label{eq:ordonMsig}
\lambda_1(\,\overline{M_{\sigma(1)}}\,)\leq
\lambda_1(\,\overline{M_{\sigma(2)}}\,)\leq
\dots\leq \lambda_1(\,\overline{M_{\sigma(n)}}\,)\;.
\end{equation}
Soit $i\in\llbracket 1, n\rrbracket$. Montrons que
\begin{equation}\label{eq:minsuccnormmax}
\lambda_i(\,\overline{M'}\,)=\lambda_1(\,\overline{M_{\sigma(i)}}\,)\;.
\end{equation}
En effet, puisque l'application $m\mapsto \lambda_1
(\,\overline{M_{\sigma(m)}} \,)$ est croissante sur $\llbracket 1,
n\rrbracket$ et comme le \mbox{$A$-module} $M_1\oplus\dots\oplus M_i$
est de rang $i$, nous avons $\lambda_i (\,\overline{M'}\,)\leq
\lambda_1 (\,\overline{M_{\sigma(i)}}\,)$.  Réciproquement, notons
$\pi_\ell:M'\ra M_{\sigma(\ell)}$ la projection canonique. La
projection canonique somme
\[
\oplus_{\ell=i}^n\pi_\ell:
M'\ra M_{\sigma(i)}\oplus\dots\oplus M_{\sigma(n)}
\]
est surjective, de noyau de rang $i-1$. Donc si $x_1,\dots, x_i\in M'$
sont $K$-linéairement indépendants, alors il existe $j\in\llbracket 1,
i\rrbracket$ et $\ell\in\llbracket i, n\rrbracket$ tels que
$\pi_\ell(x_j)\neq 0$. Puisque la norme de $\overline{M'}$ est la
norme du maximum et comme $m\mapsto \lambda_1(\,\overline{M_{\sigma(m)}}
\,)$ est croissante, nous avons donc
\[
\log_q\|x_j\|'=\max_{1\leq m\leq n}\log_q\|\pi_m(x_j)\|_{\sigma(m)}\geq
\log_q\|\pi_\ell(x_j)\|_{\sigma(\ell)}\geq \lambda_1(\,\overline{M_{\sigma(\ell)}}\,)
\geq\lambda_1(\,\overline{M_{\sigma(i)}}\,)\;.
\]
En prenant la borne inférieure sur les tels $i$-uplets $(x_1,\dots,
x_i)$, ceci montre que nous avons $\lambda_i (\,\overline{M'}\,)\geq
\lambda_1 (\,\overline{M_{\sigma(i)}}\,)$ et la formule
\eqref{eq:minsuccnormmax} en découle.

Maintenant, par le lemme \ref{lem:controlinverselong} que nous
appliquons au morphisme de $\O_{{\bf C},P_\infty}$-modules $\psi=
\varphi_{P_\infty}: (L_1)_{P_\infty} \oplus \dots\oplus
(L_n)_{P_\infty} \ra E_{P_\infty}$ induit par $\varphi$ en $P_\infty$,
et par le théorème \ref{theo:quasiscind}, l'extension $\wh K$-linéaire
$\phi_{\wh K}: M'_{\wh K}\ra M_{\wh K}$ de $\phi$ vérifie
\begin{equation}\label{eq:contorlnormop}
\interleave\,\phi_{\wh K}\interleave\leq 1\quad\text{et}\quad
\interleave(\phi_{\wh K})^{-1}\interleave\leq
q^{\operatorname{lg}\,(E_{P_\infty}/\operatorname{im}(\varphi_{P_\infty}))}\leq
q^{\frac{n(n-1)}{2\,\fff}\,(4\ggg+3\fff-3)}\;.
\end{equation}
Respectivement par l'inégalité de gauche dans la formule
\eqref{eq:contorlnormop} et la formule \eqref{eq:comparnormsuccmin},
par la formule \eqref{eq:minsuccnormmax}, par le lemme
\ref{lem:relatpentmaxsys} puisque $E^{\overline{M}_i}=L_i$ et $L_i$
est de rang $1$ donc nous avons $\mu_{\max}(L_i)=\mu_{1}(L_i)=
\deg(L_i)$, et, par le corollaire \ref{coro:pentdegdroit}, nous avons
\begin{align*}
  \lambda_i(\,\overline{M}\,)&\leq  \lambda_i(\,\overline{M'}\,)=
  \lambda_1(\,\overline{M_{\sigma(i)}}\,)\leq 
  \Big\lfloor\frac{\ggg-\mu_{\max}(L_{\sigma(i)})-1}{\fff} \Big\rfloor+1
  =\Big\lfloor\frac{\ggg-\deg(L_{\sigma(i)})-1}{\fff} \Big\rfloor+1\\ &
  \leq \Big\lfloor\frac{-\,\mu_{\sigma(i)}(E)\,+\,(n-1)(6\ggg+5\fff-5)
  \,+\,\ggg-1}{\fff}\Big\rfloor+1\;.
\end{align*}
Réciproquement, en posant $c_{\ggg,\fff,n}=\big\lfloor \frac{n(n-1)}
{2\,\fff} \,(4\ggg+3\fff-3)\big\rfloor$, par l'inégalité de droite
dans la formule \eqref{eq:contorlnormop} et la formule
\eqref{eq:comparnormsuccmin}, par la formule
\eqref{eq:minsuccnormmax}, par le lemme \ref{lem:relatpentmaxsys}
puisque $E^{\overline{M}_i}=L_i$ et par le corollaire
\ref{coro:pentdegdroit}, nous avons
\begin{align*}
  \lambda_i(\,\overline{M}\,)&\geq
  \lambda_i(\,\overline{M'}\,)-c_{\ggg,\fff,n}=
    \lambda_1(\,\overline{M_{\sigma(i)}}\,)-c_{\ggg,\fff,n}\geq
    \Big\lceil\frac{-\mu_{\max}(L_{\sigma(i)})}{\fff}
    \Big\rceil-c_{\ggg,\fff,n}\\ &
    =\Big\lceil\frac{-\deg(L_{\sigma(i)})}{\fff}
    \Big\rceil-c_{\ggg,\fff,n} \geq
\Big\lceil\frac{-\,\mu_{\sigma(i)}(E)\,-\,(n-1)(6\ggg+5\fff-5)
        \,}{\fff}\Big\rceil-c_{\ggg,\fff,n}\;.
\end{align*}
Le théorème \ref{theo:relatpenminsucc} en découle.
\cqfd

\bigskip\noindent{\bf Démonstration du théorème
  \ref{theo:reducedbasis}. } Soient $n,V,\|\;\|,\Lambda$ les termes
qui apparaissent dans l'énoncé du théorème \ref{theo:reducedbasis}, de
sorte que nous ayons $\Lambda \otimes_A \wh K=V$. Conservons les
notations $\overline{M_1}=(M_1,\|\;\|_1),\dots$, $\overline{M_n}=
(M_n,\|\;\|_n)$ et $\overline{M'}=(M',\|\;\|')$, ainsi que la
permutation $\sigma$, de la démonstration ci-dessus appliquée au
$A$-réseau normé $\overline{M}=(\Lambda,\|\;\|)$. Posons $\Lambda_1=
M_{\sigma(1)}, \dots, \Lambda_n=M_{\sigma(n)}$.  Alors $\Lambda_1,
\dots, \Lambda_n$ sont des sous-$A$-modules de rang $1$ de $\Lambda$
et nous avons $\Lambda= \Lambda_1 \oplus\dots \oplus\Lambda_n$. Posons
$c_{\ggg,\fff,n}= \big\lfloor \frac{n(n-1)} {2\,\fff} \,(4\ggg+3\fff-3)
\big\rfloor$.  Respectivement par la définition de $\overline{M}$, par
l'inégalité de gauche dans la formule \eqref{eq:contorlnormop} et la
formule \eqref{eq:comparnormsuccmin}, par la formule
\eqref{eq:minsuccnormmax}, puisque la norme de $\overline{M'}$ est la
norme maximum des normes des $\overline{M_i}$, et enfin par
l'inégalité de droite dans la formule \eqref{eq:contorlnormop} et la
formule \eqref{eq:comparnormsuccmin}, nous avons
\begin{align*}
\lambda_i(\Lambda,\|\;\|)&=\lambda_i(\,\overline{M}\,)
\leq  \lambda_i(\,\overline{M'}\,)= \lambda_1(\,\overline{M_{\sigma(i)}}\,)
=\lambda_1(\Lambda_i,\|\;\|'_{\mid \wh{K}\Lambda_i})\\&\leq 
\lambda_1(\Lambda_i,\|\;\|_{\mid \wh{K}\Lambda_i})+c_{\ggg,\fff,n}\;.
\end{align*}
Réciproquement, nous avons de même
\begin{align*}
\lambda_i(\Lambda,\|\;\|)&=\lambda_i(\,\overline{M}\,)
\geq \lambda_i(\,\overline{M'}\,)-c_{\ggg,\fff,n}=
\lambda_1(\,\overline{M_{\sigma(i)}}\,)-c_{\ggg,\fff,n}
=\lambda_1(\Lambda_i,\|\;\|'_{\mid \wh{K}\Lambda_i})- c_{\ggg,\fff,n}\\&\geq 
\lambda_1(\Lambda_i,\|\;\|_{\mid \wh{K}\Lambda_i})- c_{\ggg,\fff,n}\;.
\end{align*}
Donc nous avons $|\,\lambda_i(\Lambda,\|\;\|)- \lambda_1(\Lambda_i,
\|\;\|_{\mid \wh{K}\Lambda_i})\,|\leq c_{\ggg,\fff,n}$. Par la formule
\eqref{eq:ordonMsig}, nous avons $\lambda_1(\Lambda_1,\|\;\|_{\mid \wh
  K\Lambda_1})\leq \dots\leq \lambda_1(\Lambda_n,\|\;\|_{\mid \wh
  K\Lambda_n})$. Ceci termine la démonstration du théorème
\ref{theo:reducedbasis}.  \cqfd

{\small \bibliography{../biblio}}

\begin{thebibliography}{BAPP}

\bibitem[And]{Andre09}
Y.~André.
\newblock {\it Slope filtrations}.
\newblock {Confluentes Math. {\bf 1} (2009) 1--85}.

\bibitem[Art]{Arthur78}
J.~G.~Arthur.
\newblock {\it A trace formula for reductive groups. I. Terms associated
to classes in $G({\bf Q})$}.
\newblock {Duke Math. J. {\bf 45} (1978) 911–952}.

\bibitem[BKLP]{BanKimLimPau24}
  G.~Bang, T.~Kim,  S.~Lim and F.~Paulin.
\newblock {\it Parametric geometry of number over general function fields}.
\newblock {In preparation}.
 
\bibitem[Bos1]{Bost96}
J.~B.~Bost.
\newblock {\it Périodes et isogénies des variétés abéliennes
sur les corps de nombres}.
\newblock {Astérisque {\bf 237}  (1996), Séminaire Bourbaki Vol. 1994/1995.
Exp. No 795,  115--161}.
  
\bibitem[Bos2]{Bost20}
J.~B.~Bost.
\newblock {\it  Réseaux euclidiens, séries thêta et pentes [d'après
    W. Banaszczyk, O. Regev, D. Dadush, N. Stephens-Davidowitz, …]}.
\newblock {Astérisque {\bf 422} (2020), Séminaire Bourbaki Vol.
  2018/2019.  Exp. No. 1151, 1--59}.

\bibitem[BosC]{BosCha24}
J.~B.~Bost and F.~Charles.
\newblock {\it  Infinite Dimensional Geometry of Numbers:
    Hermitian Quasi-coherent Sheaves and Theta Finiteness}.
\newblock {Book preprint, available on
    https\!:// www.math.ens.psl.eu/$\sim$charles/inf$\_$dim$\_$1}.

\bibitem[BdS]{BreSax23}
E.~Breuillard, N.~de Saxcé.
\newblock {\it  A subspace theorem for manifolds}.
\newblock {J. Eur. Math. Soc. {\bf 26} (2023) 4273--4313}.

\bibitem[BPP]{BroParPau19}
A.~Broise-Alamichel, J.~Parkkonen and F.~Paulin.
\newblock {\it Equidistribution and counting under equilibrium states in
  negative curvature and trees. Applications to non-Archimedean Diophantine
  approximation}.
\newblock {With an Appendix by J.~Buzzi. Prog. Math. {\bf 329}, Birkhäuser,
  2019}.

\bibitem[Che1]{Chen10a}
Huayi Chen.
\newblock {\it Convergence des polygones de Harder-Narasimhan}.
\newblock {Mémoires  Soc. Math. France {\bf 120} (2010)}.

\bibitem[Che2]{Chen10b}
Huayi Chen.
\newblock {\it Harder-Narasimhan categories}.
\newblock {J. Pure Appl. Algebra {\bf 214} (2010) 187--200}.

\bibitem[CM]{CheMor20}
Huayi Chen and A.~Moriwaki.
\newblock {\it Arakelov geometry over adelic curves}.
\newblock {Lecture Notes in Math. {\bf 2258}, Springer,  2020}.

\bibitem[Cor]{Cornut18}
C.~Cornut.
\newblock {\it On Harder–Narasimhan filtrations and their
  compatibility with tensor products}.
\newblock {Confluentes Math. {\bf 10} (2018) 3--49}.

\bibitem[deS]{deSaxce23}
N.~de Saxcé.
\newblock {\it Non-divergence in the space of lattices}.
\newblock {Groups Geom. Dyn. {\bf 17} (2023) 993--1003}.

\bibitem[FalW]{FalWus94}
G.~Faltings and G.~Wüstholz.
\newblock {\it Diophantine approximations on projective spaces}.
\newblock {Inventiones Mathematicae {\bf 116} (1994) 109--138}.

\bibitem[Far]{Fargues10}
L.~Fargues.
\newblock {\it La filtration de Harder-Narasimhan des schémas
  en groupes finis et plats}.
\newblock {J. reine angew. Math. [Crelle's Journal] {\bf 645} (2010) 1--39}.

\bibitem[Gau]{Gaudron08}
É.~Gaudron.
\newblock {\it Pentes des fibrés vectoriels adéliques sur un corps global}.
\newblock {Rend. Semin. Mat. Univ. Padova {\bf 119} (2008) 21--95}.

\bibitem[GI]{GolIwa63}
O.~Goldman and N.~Iwahori.
\newblock {\it The space of p-adic norms}.
\newblock {Acta Math. {\bf 109} (1963) 137--177}.

\bibitem[Gos]{Goss96}
D.~Goss.
\newblock {\it Basic structures of function field arithmetic}.
\newblock {Erg. Math. Grenz. {\bf 35}, Springer Verlag 1996}.

\bibitem[Gra]{Grayson84}
D.~Grayson.
\newblock {\it Reduction theory using semistability}.
\newblock {Comm. Math. Helv {\bf 59} (1984) 600--634}.

\bibitem[Gri]{Grieve23}
N.~Grieve.
\newblock {\it Vertices of the Harder and Narasimhan polygons
  and the laws of large numbers}.
\newblock {Can. Math. Bull. {\bf 66} (2023) 340--357}.

\bibitem[Gro]{Grothendieck57}
A.~Grothendieck.
\newblock {\it Sur la classification des fibrés holomorphes
  sur la sphère de Riemann}.
\newblock {Amer. J. Math. {\bf 79} (1957) 121--138}.

\bibitem[Har]{Harder69}
G.~Harder.
\newblock {\it Minkowskische Reduktionstheorie über Funktionenkörpern}.
\newblock {Invent. Math. {\bf 7} (1969) 33--54}.

\bibitem[HN]{HarNar75}
G.~Harder and M.~S.~Narasimhan.
\newblock {\it On the Cohomology Groups of Moduli Spaces
of Vector Bundles on Curves}.
\newblock {Math. Ann {\bf 212} (1975) 215--248}.

\bibitem[HuL]{HuyLeh10}
  D.~Huybrechts and M.~Lehn.
\newblock {\it The Geometry of Moduli Spaces of Sheaves}.
\newblock {2nd ed. Cambridge Univ. Press 2010}.

\bibitem[KLP]{KimLimPau23}
T.~Kim, S.~Lim and F.~Paulin.
\newblock {\it On Hausdorff dimension in inhomogeneous Diophantine
  approximation over global function fields}.
\newblock {J. Numb. Theo. {\bf 251} (2023) 102--146}.

\bibitem[KlST]{KleShiTom17}
D.~Kleinbock, R.~Shi and G.~Tomanov.
\newblock {\it S-adic version of Minkowski's geometry of numbers
  and Mahler's compactness criterion}.
\newblock {J. Numb. Theo. {\bf 174} (2017) 150--163}.

\bibitem[LaLS]{LagLenSch90}
J.~Lagarias, H.~Lenstra, Jr. and C.-P.~Schnorr.
\newblock {\it Korkin-Zolotarev bases and successive minima of
  a lattice and its reciprocal lattice}.
\newblock {Combinatorica {\bf 10} (1990) 333--348}.

\bibitem[Len]{Lenstra85}
A.~K.~Lenstra.
\newblock {\it Factoring multivariate polynomials over finite fields}.
\newblock {J. Comput. System Sci. {\bf 30} (1985) 235--248}.

\bibitem[Li]{Li24}
Yao Li.
\newblock {\it Categorification of Harder–Narasimhan theory
  via slope functions valued in totally ordered sets}.
\newblock {Manuscripta Math. {\bf 173} (2024) 1233--1271}.

\bibitem[Mah]{Mahler41}
K.~Mahler.
\newblock {\it An analogue to Minkowski's geometry of numbers in a field
         of series}.
\newblock {Ann. of Math. {\bf 42} (1941) 488--522}.

\bibitem[Par]{Parreau00}
A.~Parreau.
\newblock {\it  Immeubles affines: construction par les normes et
  étude des isométries}.
\newblock {Contemp. Math. {\bf 262},
Amer. Math. Soc. 2000, 263--302}.

\bibitem[PR]{PoeRoy23}
A.~Poëls and D.~Roy.
\newblock {\it Parametric geometry of numbers over a number field
and extension of scalars}.
\newblock {Bull. Soc. Math. France {\bf 151} (2023) 257--303}.

\bibitem[Ros]{Rosen02}
M.~Rosen.
\newblock {\it Number theory in function fields}.
\newblock {Grad. Texts Math. {\bf 210}, Springer Verlag, 2002}.

\bibitem[RW]{RoyWal17}
D.~Roy and M.~Waldschmidt.
\newblock {\it Parametric geometry of numbers in function fields}.
\newblock {Mathematika {\bf 63} (2017) 1114--1135}.

\bibitem[Ser1]{Serre55}
J.-P.~Serre.
\newblock {\it Faisceaux Algebriques Coherents}.
\newblock {Ann. of Math. {\bf 61} (1955) 197--2783}.

\bibitem[Ser2]{Serre83}
J.-P.~Serre.
\newblock {\it Arbres, amalgames, SL$_2$}.
\newblock {3\`eme \'ed. corr., Ast\'erisque {\bf 46}, Soc. Math. France, 1983}.

\bibitem[Stu]{Stuhler76}
U.~Stuhler.
\newblock {\it Eine Bemerkung zur Reduktionstheorie quadratischer Formen}.
\newblock {Arch. Math. {\bf 27} (1976) 604--610}.

\bibitem[Thu]{Thunder94}
J.~L.~Thunder.
\newblock {\it Siegel's lemma for function fields}.
\newblock {Michigan Math. J. {\bf 42}  (1995) 147--162}.

\bibitem[Wei]{Weil95}
A.~Weil.
\newblock {\it Basic number theory}.
\newblock {Classics in Math, Springer Verlag, 1995}.

\end{thebibliography}

{\small \bigskip\noindent \begin{tabular}{l} 
  Laboratoire de math\'ematique d'Orsay,
  UMR 8628  CNRS\\ Universit\'e Paris-Saclay, 91405
  ORSAY Cedex, FRANCE
  \\ {\it jean-benoit.bost@universite-paris-saclay.fr}
  \\{\it frederic.paulin@universite-paris-saclay.fr}
\end{tabular}
}

\end{document}